\documentclass[10pt]{article}

\usepackage{amsmath,amsfonts,graphics,graphicx,bm,enumerate,epsfig,psfrag,subfig,longtable}
\usepackage{palatino,cite}
\usepackage{color,soul} 
\usepackage{xcolor}
\usepackage[mathscr]{euscript} 
\usepackage{mathabx} 
\usepackage{float}
\usepackage[colorlinks=true, pdfstartview=FitV, linkcolor=black, citecolor= black, urlcolor= black]{hyperref}
\usepackage{footnpag}			       

\newcommand{\m}[1]{{\bf{#1}}}

\newcommand{\tr}{^{\sf T}}
\newcommand{\C}[1]{{\cal {#1}}}

\newcommand{\w}[1]{\omega_{#1}}
\newcommand{\e}[1]{\epsilon_{#1}}

\newcommand{\edot}[1]{\dot{\epsilon}_{#1}}

\title{\bf{End-to-End Ascent-Entry Mission Performance Optimization Using Gaussian Quadrature Collocation}}

\author{Alexander T.~Miller\footnote{Ph.D.~Student, NDSEG Fellow, Department of Mechanical and Aerospace Engineering.  E-mail:  alexandertmiller@ufl.edu.} \\ Anil~V.~Rao\footnote{Professor, Erich Farber Faculty Fellow, and University Term Professor, Department of Mechanical and Aerospace Engineering.  E-mail:  anilvrao@ufl.edu.  Associate Fellow AIAA.  Corresponding Author.} \\\\ {\em University of  Florida} \\ {\em Gainesville, FL 32611}}

\begin{document}

\date{}
\maketitle{}
\renewcommand{\baselinestretch}{1}\normalsize\normalfont

 \begin{abstract}
   The performance optimization for a combined ascent-entry mission subject to constraints on heating rate and heating load is studied.  The ascent vehicle is modeled as a three-stage rocket that places the vehicle onto a suborbital exo-atmopheric trajectory after which the vehicle undergoes an unpowered entry and descent to a vertically downward terminal condition.  The entry vehicle is modeled as a high lift-to-drag ratio vehicle that is capable of withstanding high levels of thermal and structural loads.  A performance index is designed to improve control margin while attenuating phugoid oscillations during atmospheric entry.  Furthermore, a mission corresponding to a prototype launch and target point is used in this study.  The trajectory optimization problem is formulated as a multiple-phase optimal control problem, and the optimal control problem is solved using an adaptive Gaussian quadrature collocation method.  A key aspect of the optimized trajectories is that, for particular ranges of maximum allowable heating rate and heating load during entry, relatively small adjustments made during ascent can potentially decrease the control effort required during atmospheric entry.  Outside of these ranges for maximum allowable heating rate and heating load, however, it is found that the required control effort increases and eventually saturates the commanded angle of attack upon initial descent.  The key features of the optimized trajectories and controls are identified, and the approach developed in this paper provides a systematic method for end-to-end ascent-entry trajectory optimization.  
\end{abstract}

\section*{Nomenclature}
\renewcommand{\baselinestretch}{1}\normalsize\normalfont
\begin{longtable}{lcl}
$C$        & $=$ & penalty term constant\\
$C_D$        & $=$ & drag coefficient\\
$C_{D0}$   & $=$ & zero-lift drag coefficient\\
$C_L$        & $=$ & lift coefficient\\
$D$        & $=$ & drag force magnitude, kN\\
$g_0$     & $=$ & standard acceleration due to gravity, km/s$^2$\\
$h$        & $=$ & altitude over spherical Earth, km\\
$h_{\textrm{atm}}$ & $=$ & pierce point altitude, km\\
$h_{\textrm{peak},\max}$ & $=$ & maximum allowable peak altitude, km\\
$h_{\textrm{peak},\min}$ & $=$ & minimum allowable peak altitude, km\\
$I_{SP}$ & $=$ & specific impulse, s\\
$\C{J}$  & $=$ & cost\\
$\C{J}^{(p)}$  & $=$ & phase $p$ cost\\
$k$        & $=$ & penalty term design variable\\
$K$        & $=$ & drag polar parameter\\
$L$        & $=$ & lift force magnitude, kN\\
$\C{L}$ & $=$ & penalty term\\
$m$        & $=$ & mass, kg\\
$m_{fairing}$ & $=$ & fairing mass, kg\\
$m_{S2}$ & $=$ & mass at stage 2 ignition, kg\\
$m_{S3}$ & $=$ & mass at stage 3 ignition, kg\\
$M_1$       & $=$ & Mission 1\\
$M_2$       & $=$ & Mission 2\\
$M_3$       & $=$ & Mission 3\\
$n$        & $=$ & sensed acceleration, g\\
$n_{\max}$ & $=$ & maximum allowable sensed acceleration, g\\
$q$        & $=$ & dynamic pressure, kPa\\
$q_{\max}$ & $=$ & maximum allowable dynamic pressure, kPa\\
$q_{\min}$ & $=$ & minimum allowable dynamic pressure, kPa\\
$Q$        & $=$ & heating load, MJ/m$^2$ \\
$Q_{\max}$ & $=$ & maximum allowable heating load, MJ/m$^2$ \\
$\dot{Q}$ & $=$ & heating rate, MW/m$^2$\\
$\dot{Q}_{\max}$ & $=$ & maximum allowable heating rate, MW/m$^2$\\
$R_e$     & $=$ & geocentric radius, km\\
$S$        & $=$ & reference area, m$^2$\\
$t$         & $=$ & time, s\\
$t_0^{(p)}$  & $=$ & initial time of phase $p$\\
$t_f^{(p)}$  & $=$ & terminal time of phase $p$\\
$t_{\textrm{fairing}}$ & $=$ & fairing separation time, s\\
$t_{S1}$ & $=$ & stage 1 burnout time, s\\
$t_{S2}$ & $=$ & stage 2 burnout time, s\\
$t_{S3}$ & $=$ & stage 3 burnout time, s\\
$T$        & $=$ & thrust, kN\\
$u_1$      & $=$ & slack variable, rad/s\\
$u_2$     & $=$ & slack variable, rad/s\\
$u_{\alpha}$ & $=$ & angle of attack rate, deg/s or rad/s\\
$u_{\alpha,\max}$ & $=$ & maximum allowable angle of attack rate, deg/s or rad/s\\
$u_{\sigma}$ & $=$ & bank angle rate, deg/s or rad/s\\
$u_{\sigma,\max}$ & $=$ & maximum allowable bank angle rate, deg/s or rad/s\\
$v$        & $=$ & Earth-relative speed, km/s\\
$v_c$     & $=$ & Earth radius circular speed, km/s\\
$\m{y}$ & $=$ & state vector\\
$\alpha$ & $=$ & angle of attack, deg or rad\\
$\bar{\alpha}$ & $=$ & angle of attack setpoint, deg or rad\\
$\alpha_{\max}$ & $=$ & maximum allowable angle of attack, deg or rad\\
$\gamma$ & $=$ & flight-path angle, deg or rad\\
$\e{1}$ & $=$ & Euler parameter\\
$\e{2}$ & $=$ & Euler parameter\\
$\e{3}$ & $=$ & Euler parameter\\
$\eta$  & $=$ & Euler parameter\\
$\theta$ & $=$ & latitude, deg or rad\\
$\kappa$ & $=$ & heating rate parameter, MW/m$^2$ \\
$\mu_e$ & $=$ & Earth gravitational parameter, km$^3$/s$^2$\\
$\rho$     & $=$ & atmospheric density, kg/m$^3$\\
$\rho_0$  & $=$ & atmospheric density at sea level, kg/m$^3$\\
$\sigma$ & $=$ & bank angle, deg or rad\\
$\phi$     & $=$ & longitude, deg or rad\\
$\psi$     & $=$ & azimuth angle, deg or rad\\
$\w{e}$  & $=$ & Earth rotation rate, deg/s or rad/s\\
$\w{1}$  & $=$ & angular velocity component, deg/s or rad/s\\
$\w{2}$  & $=$ & angular velocity component, deg/s or rad/s\\
$\w{3}$  & $=$ & angular velocity component, deg/s or rad/s
\end{longtable}
\renewcommand{\baselinestretch}{2}\normalsize\normalfont

\section{Introduction}
The strategic importance of highly maneuverable, high-speed aircraft cannot be understated.  High-speed vehicles described in Ref.~\cite{Richie1999} enable rapid and precise global reach for a plethora of missions ranging from support operations to target strikes.  Such missions pose challenging vehicle design, trajectory optimization, and guidance problems due to the tight constraints experienced during entry (limits on heat rate, heat load, g-load, dynamic pressure, etc.).  Thus, the past two decades have seen a wide array of research into constrained entry trajectory optimization and guidance for high-speed glide vehicles (for instance, Refs.~\cite{Lu2003, Lu2014, Lu2016, Clarke1, Clarke2, Jorris1, Jorris2, Jorris3, Tian2011, Zhao2013, Zhao2014, Moshman2014, Dong2017}).  The introduction of Ref.~\cite{Lu2014} and the work of Ref.~\cite{ Zhao2014} provides a summary of the recent work in this area.

Compared with the studies performed on high-speed entry vehicles, less research has focused on trajectory optimization for the combined ascent-entry problem \cite{Li2009, Lu2013, Rizvi1, Rizvi2, Rizvi3, Miller1}.  Reference~\cite{Li2009} employs a direct shooting method to generate a maximum downrange trajectory for a high-speed vehicle boosted from a three-stage rocket.  A predictor-corrector method is employed in Ref.~\cite{Lu2013} for gliding guidance of high lift-to-drag ratio hypersonic vehicles launched from a multi-stage rocket.  References~\cite{Rizvi1, Rizvi2, Rizvi3} utilize an adaptive pseudospectral method to explore trade-offs in maximum downrange or crossrange performance for a variety of entry vehicles boosted from a two-stage rocket.  Finally, Ref.~\cite{Miller1} employs a model of the space shuttle to study an ascent-entry problem that primarily seeks to maximize the payload mass delivered to low Earth orbit (LEO).  Taken together, the aforementioned works point to the possibility of end-to-end mission planning for a wide variety of ascent-entry missions.  In addition, the previous works indicate the critical role ascent trajectory shaping plays in improving performance during entry due in the presence of key constraints (heating rate, sensed acceleration load, etc.) placed on the vehicle during entry.

Motivated by the previous research, this paper studies end-to-end mission planning for a hypothetical mission involving a high-speed vehicle with performance characteristics similar to that given in Ref.~\cite{Richie1999}.  The mission begins with the launch of a three-stage rocket to boost the vehicle to suborbital speeds.  The mission ends with the vehicle vertically striking the target at the desired endpoint conditions.  Performance is graded by an objective that seeks to maintain wide control margins at all times, and variations in performance are analyzed as the heating rate and heating load requirements during entry are tightened.  Furthermore, the addition of a penalty term to the Lagrange cost is studied in an attempt to reduce phugoid oscillations during entry.  Finally, all optimal trajectories are obtained via Legendre-Guass-Radau (LGR) collocation \cite{Garg1, Darby3}, chosen for its desirable convergence properties \cite{HagerHouRao15c, HagerHouRao16a} in addition to its ability to handle the large number of constraints, highly nonlinear dynamics, and integrals involved in the problem.

This research takes a different approach from the previous research into combined ascent-entry trajectory optimization.  First, a Gaussian quadrature method is used (in contrast to Refs.~\cite{Li2009, Lu2013}).  Gaussian quadrature methods have grown in popularity in recent years because they allow for a highly general problem formulation, path constraints and integral constraints are enforced with ease, and high accuracy solutions can be obtained even when the initial guess is poor.  Next, this paper studies a hypothetical ground strike mission whereas previous work has largely focused on maximum downrange/crossrange trajectories (see Refs.~\cite{Li2009, Rizvi1, Rizvi2, Rizvi3}.  As such, the endpoint conditions and objective functional are fundamentally different from those explored in Refs.~\cite{Li2009, Rizvi1, Rizvi2, Rizvi3, Miller1}.  Of particular note, vertical flight endpoint conditions are investigated in this paper (challenging due to singularities in the equations of motion and because the bank angle becomes undefined in vertical flight).  Finally, more realistic models are employed.  The Earth is modeled as a rotating sphere, aerodynamic coefficients are functions of both angle of attack and mach number, and all common constraints (heating rate, heating load, sensed acceleration, etc.) are enforced.

Overall, this work takes a step towards applying Gaussian quadrature methods to real-time, outer-loop guidance of high-speed, combined ascent-entry strike missions.  For instance, the method of Ref.~\cite{Dennis2019} could be employed to truncate the ascent-entry optimal control problem stated in this paper and generate guidance updates for the unexpired phases of flight.  Such an approach requires that unmodeled perturbations do not push the vehicle off course too much and cause the problem to become infeasible.  Thus, this paper's analysis of heating rate and heating load limitations as well as the study of performance indices that maximize control margin while limiting entry stresses due to phugoid oscillations will provide valuable information for future researchers.

The remainder of this paper is organized as follows.  Section~\ref{sect:phases} describes the mission and introduces each of the phases.  The physical models employed for the Earth, ascent vehicle, and entry vehicle are detailed in Section~\ref{sect:model}.  Next, the trajectory optimization problem is formed in Section~\ref{sect:ProblemFormulation} followed by several performance studies in Section~\ref{sect:results}.  Finally, Section~\ref{sect:discussion} discusses the results of the aforementioned studies and Section~\ref{sect:Conclusion} draws conclusions on this research.


\section{Ascent-Entry Mission Profile \label{sect:phases}}

\begin{figure}[hbt!]
  \centering
  \includegraphics[width=.475\textwidth]{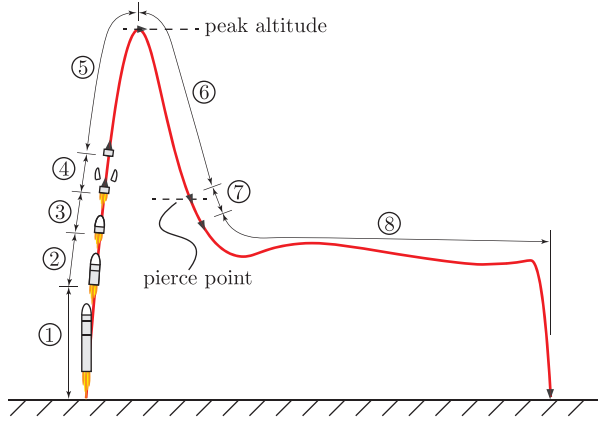}
  
  \caption{Visualization of the phase sequence for the combined ascent-entry mission.}\label{fig:phases}
\end{figure}

In this research, a combined ascent-entry mission profile is developed.  The mission consists of eight sequential phases.  Phases 1 through 4 are boost phases, where  phase 1 starts immediately after launch and phase 4 terminates at burnout of the last ascent engine.  Next, phases 5 and 6 occur in exo-atmospheric flight, where payload separation occurs in phase 5 while orientation adjustments prior to entry occur in phase 6.  Next, phases 7 and 8 comprise the atmospheric entry phases with phase 7 starting at pierce point of the atmosphere while phase 8 terminates at target impact.  Using the aforementioned descriptions, the individual phases are each summarized as follows:
\begin{enumerate}
	\item \textbf{Phase 1:}~~~Stage 1, tower cleared to engine burnout,
	\item \textbf{Phase 2:}~~~Stage 2, ignition to engine burnout,
	\item \textbf{Phase 3:}~~~Stage 3, ignition to fairing separation,
	\item \textbf{Phase 4:}~~~Stage 3, fairing separation to engine burnout,
	\item \textbf{Phase 5:}~~~Exo-atmospheric coast to peak altitude,
	\item \textbf{Phase 6:}~~~Peak altitude to pierce point,
	\item \textbf{Phase 7:}~~~Atmospheric entry (low dynamic pressure),
	\item \textbf{Phase 8:}~~~Atmospheric entry (high dynamic pressure).
\end{enumerate}
Figure~\ref{fig:phases} depicts the phase sequence.  Finally, it is noted that noting that phases 1 and 8 start and terminate, respectively, in vertical flight.

\newpage 
\section{Physical Model \label{sect:model}}
\subsection{Earth}
For this study motion the motion of a point mass is considered over a spherical rotating Earth.  Furthermore, the atmosphere is assumed to be fixed to the Earth and is modeled by the 1962 US Standard Atmosphere \cite{Sissenwine1962}.  Ambient atmospheric density and speed of sound values are calculated via interpolation of the aforementioned atmospheric data\cite{Fritsch1980}.  Table~\ref{tab:Earth-data} provides the relevant data for the modeling assumptions used in this research.  

\begin{table}[h]
  \centering
  \caption{Physical Constants for the Earth.\label{tab:Earth-data}}
  \renewcommand{\baselinestretch}{1}\normalsize\normalfont
 \begin{tabular}{lccl}\hline
    Description & Symbol & Units & Value \\\hline
    Gravitational Parameter & $\mu_{e}$ & $\textrm{km}^3/\textrm{s}^2$ & $3.986004405 \times 10^5$ \\
    Radius & $R_e$ & km & $6.378166 \times 10^3$ \\
    Rotation Rate & $\omega_e$ & $\textrm{rad}/\textrm{s}$ & $7.292115856 \times 10^{-5}$ \\
    Standard Gravity & $g_0$ & $\textrm{km}/\textrm{s}^2$ & $9.8066498 \times 10^{-3}$ \\\hline
  \end{tabular}
\end{table}

\subsection{Boost Vehicle}
The boost vehicle is modeled as a three stage rocket \cite{Min456UsersGuide, Min4DataSheet}.  All three stages employ solid rocket motors and the pertinent data for each stage is contained in Table~\ref{tab:MinIV-data}.  Note that thrust is assumed to be constant.  In addition, the fairing and payload masses are given as $400$ kg and $3000$ kg, respectively.

\begin{table}[h]
  \centering
  \caption{Boost Vehicle Data.\label{tab:MinIV-data}}
  \renewcommand{\baselinestretch}{1}\normalsize\normalfont
 \begin{tabular}{lccrrr}\hline
    Description & Symbol & Units & Stage 1 & Stage 2 & Stage 3 \\\hline
    Burn Time & - & s & 56.4 & 60.7 & 72.0 \\
    Mass (empty) & - & kg & 3630 & 3170 & 630 \\
    Mass (fuel) & - & kg & 45360 & 24500 & 7080 \\
    Mass (total) & - & kg & 48990 & 27670 & 7710 \\
    Reference Area & $S$ & $\textrm{m}^2$ & 4.307 & 4.307 & 4.307 \\
    Specific Impulse & $I_{SP}$ & s & 282 & 309 & 300 \\
    Thrust & $T$ & kN & 2224.1 & 1222.9 & 289.1 \\\hline
  \end{tabular}
\end{table}

\begin{figure}[hbt!]
  \centering
  \begin{tabular}{lr}
  \subfloat[Lift coefficient, $C_L$ vs. angle of attack, $\alpha$.]{\includegraphics[width=.475\textwidth]{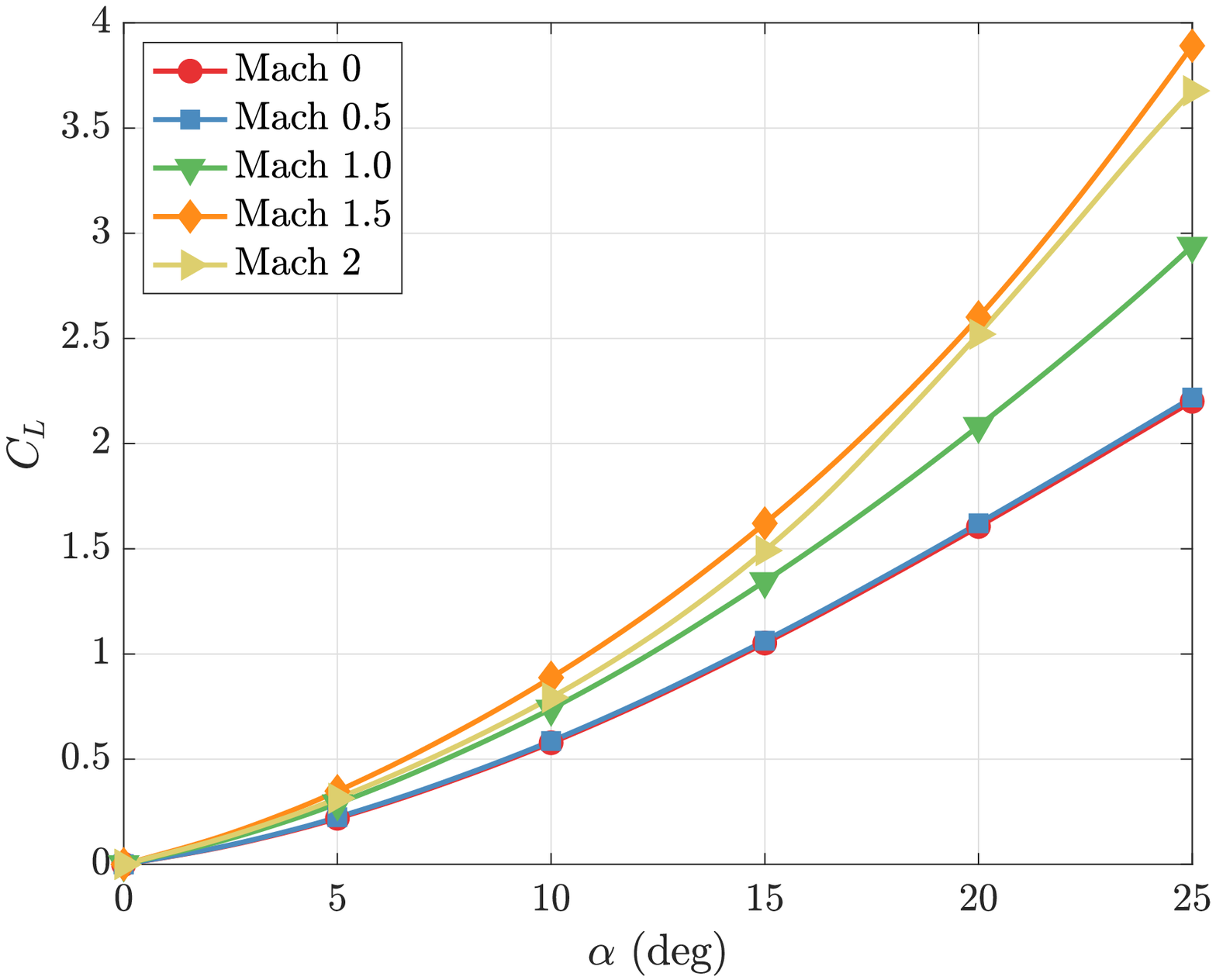}}
  &
  \subfloat[Drag coefficient, $C_D$ vs. angle of attack, $\alpha$.]{\includegraphics[width=.475\textwidth]{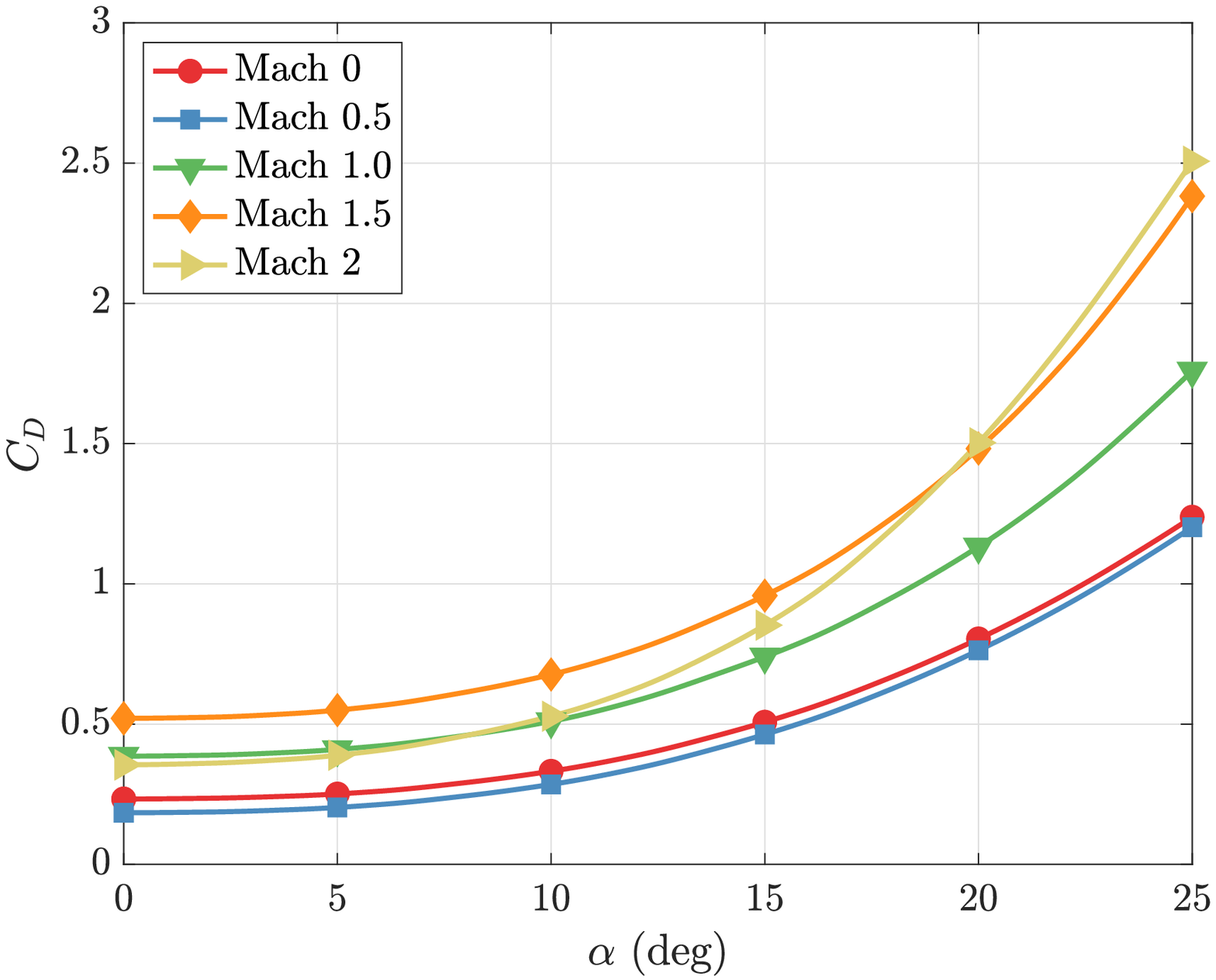}}
  \end{tabular}\\
  \subfloat[Lift-to-drag ratio, $L/D$ vs. angle of attack, $\alpha$.]{\includegraphics[width=.475\textwidth]{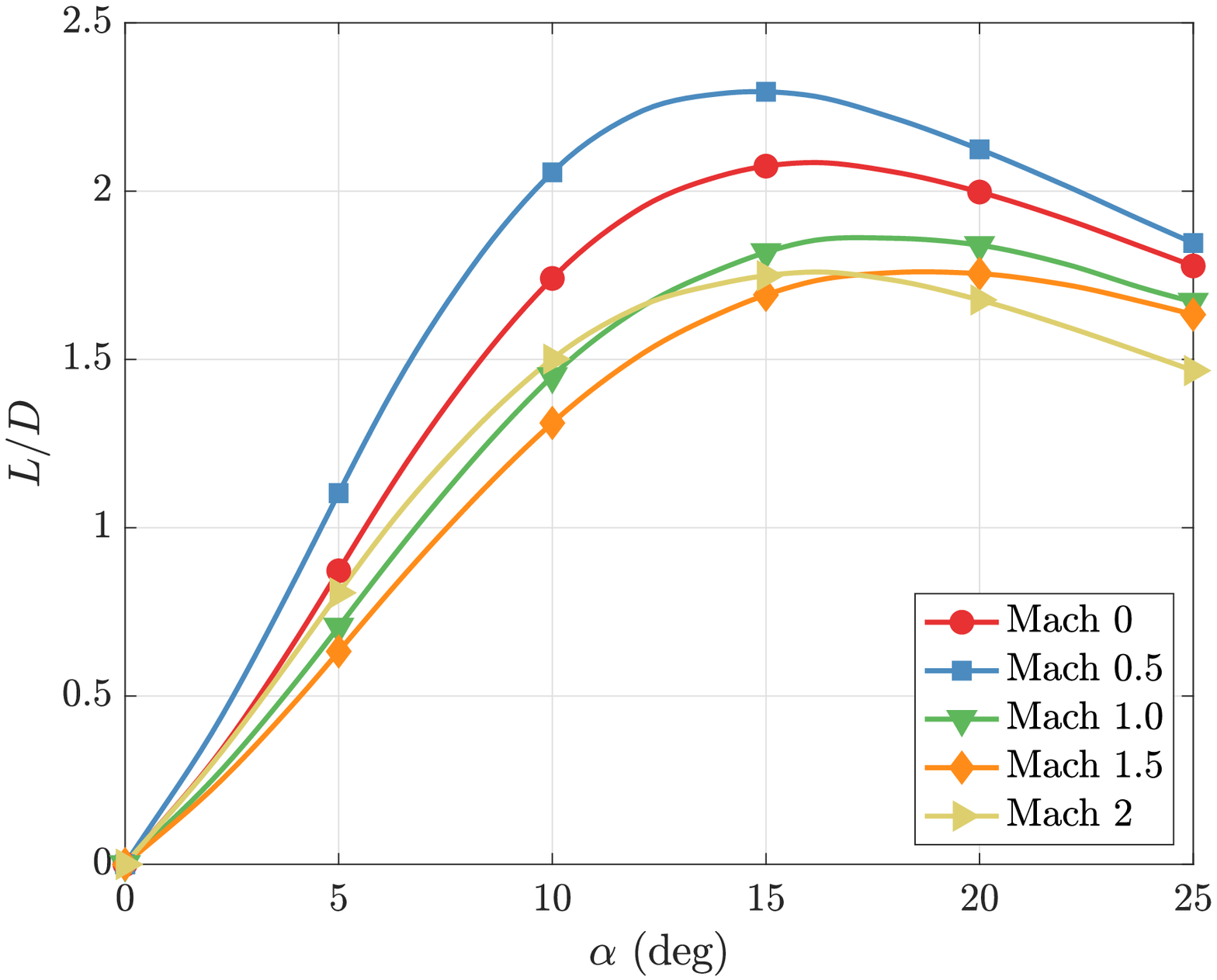}}
  
  \caption{Aerodynamic model for stage 1 at low Mach numbers.}\label{fig:MinIV-LD-low}
\end{figure}

\begin{figure}[hbt!]
  \centering
  \begin{tabular}{lr}
  \subfloat[Lift coefficient, $C_L$ vs. angle of attack, $\alpha$.]{\includegraphics[width=.475\textwidth]{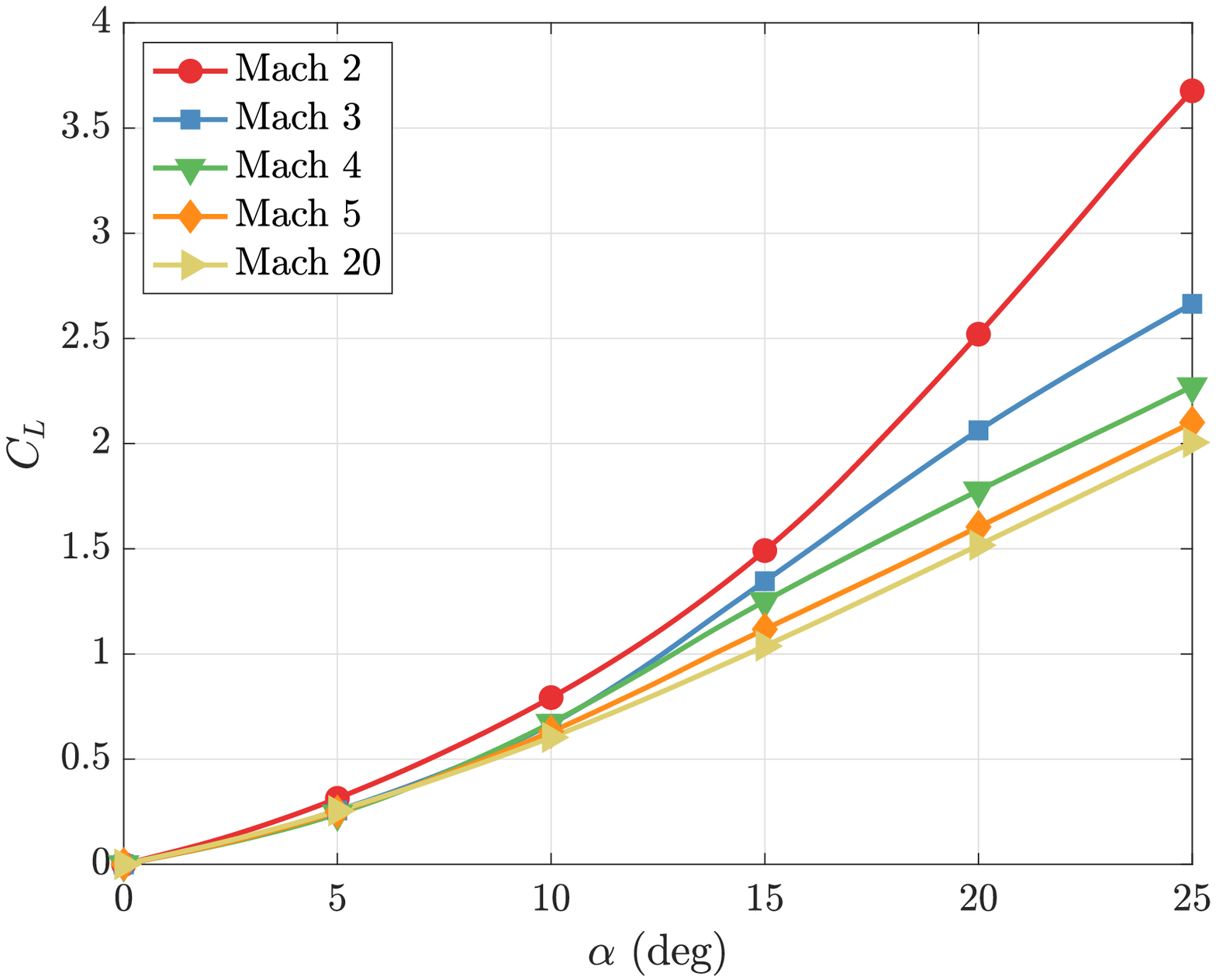}}
  &
  \subfloat[Drag coefficient, $C_D$ vs. angle of attack, $\alpha$.]{\includegraphics[width=.475\textwidth]{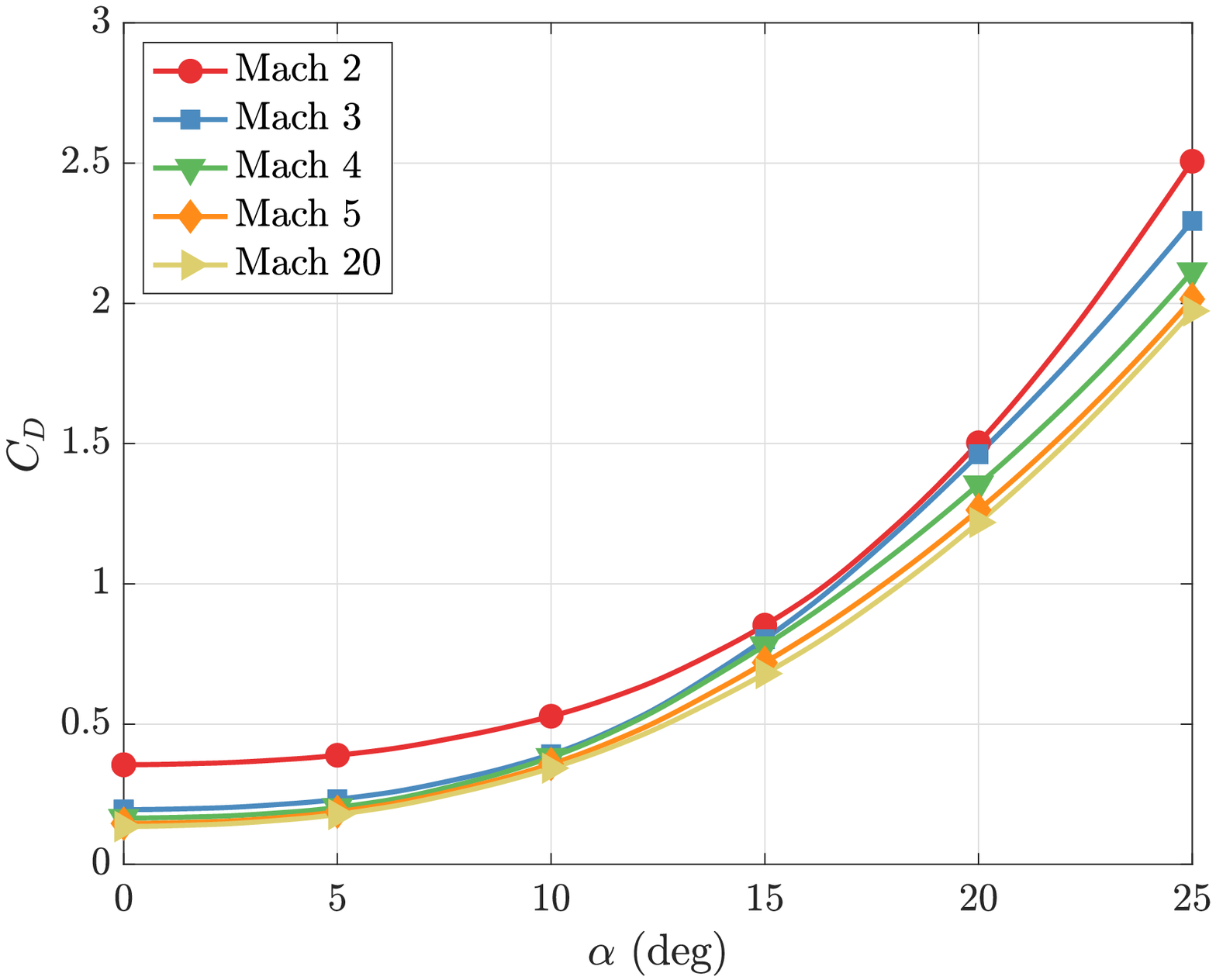}}
  \end{tabular}\\
  \subfloat[Lift-to-drag ratio, $L/D$ vs. angle of attack, $\alpha$.]{\includegraphics[width=.475\textwidth]{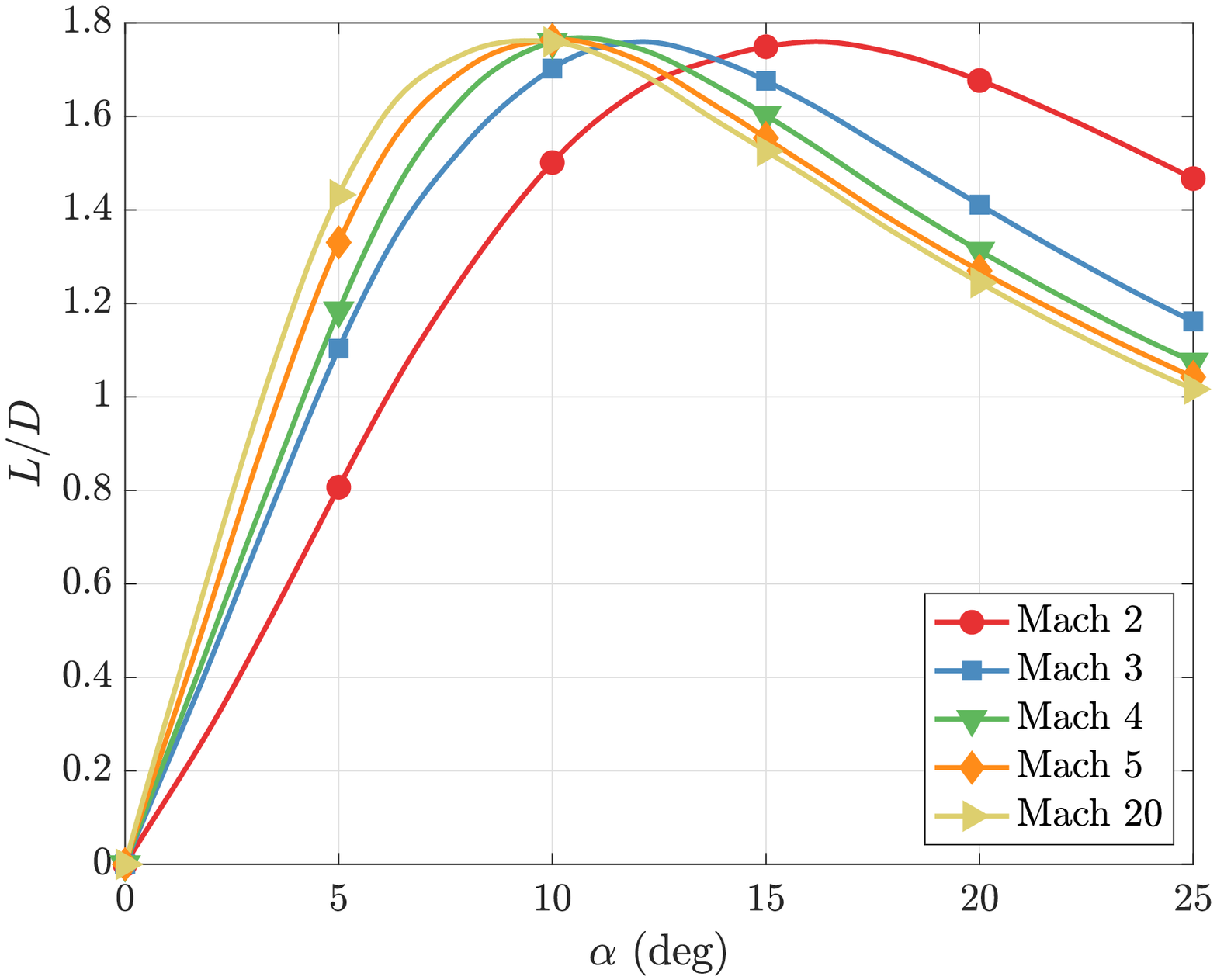}}
  
  \caption{Aerodynamic model for stage 1 at high Mach numbers.}\label{fig:MinIV-LD-high}
\end{figure}

Next, aerodynamic models are created for each stage using Digital Datcom\cite{DatcomPDAS, Finck1978, DatcomV1, DatcomV2, DatcomV3}.  Digital Datcom takes in the vehicle geometry and outputs estimates of the lift and drag coefficients, denoted $C_L$ and $C_D$ respectively, as both angle of attack, denoted $\alpha$, and Mach number are varied.  The data obtained by Digital Datcom provides an educated guess of the boost vehicle's aerodynamic characteristics, noting that the emphasis here is on qualitative behavior.  In keeping with the qualitative emphasis, it is noted that $C_L$ data generated past Mach 1.2 has been scaled down in order to maintain a reasonable lift-to-drag ratio in the supersonic and hypersonic regimes.  Specific values of $C_L$ and $C_D$ are obtained from the adjusted data via 2D interpolation\cite{makima}.  The resulting lift and drag coefficient models are displayed in Figs.~\ref{fig:MinIV-LD-low} and \ref{fig:MinIV-LD-high} for stage 1, noting that the models obtained for stages 2 and 3 are qualitatively similar.

\newpage 
\subsection{Entry Vehicle}

The model for the unpowered entry vehicle is based on the model employed in Ref.~\cite{Jorris1}.  The lift and drag coefficients are modeled as functions of the angle of attack and Mach number based on data obtained from Ref.~\cite{Jorris1} (originally from Ref.~\cite{Phillips1}).  However, the data from Ref.~\cite{Jorris1} only includes values of the angle of attack at $\alpha = \{10, 15, 20\}~\textrm{deg}$ and the angle of attack range employed in this research is $\alpha \in [0,25]~\textrm{deg}$.  Table~\ref{tab:CAV-data} provides the relevant data and vehicle aerodynamic model for the entry vehicle model used in this study.   

\begin{figure}[hbt!]
  \centering
  \begin{tabular}{lr}
  \subfloat[Lift coefficient, $C_L$ vs. angle of attack, $\alpha$.]{\includegraphics[width=.475\textwidth]{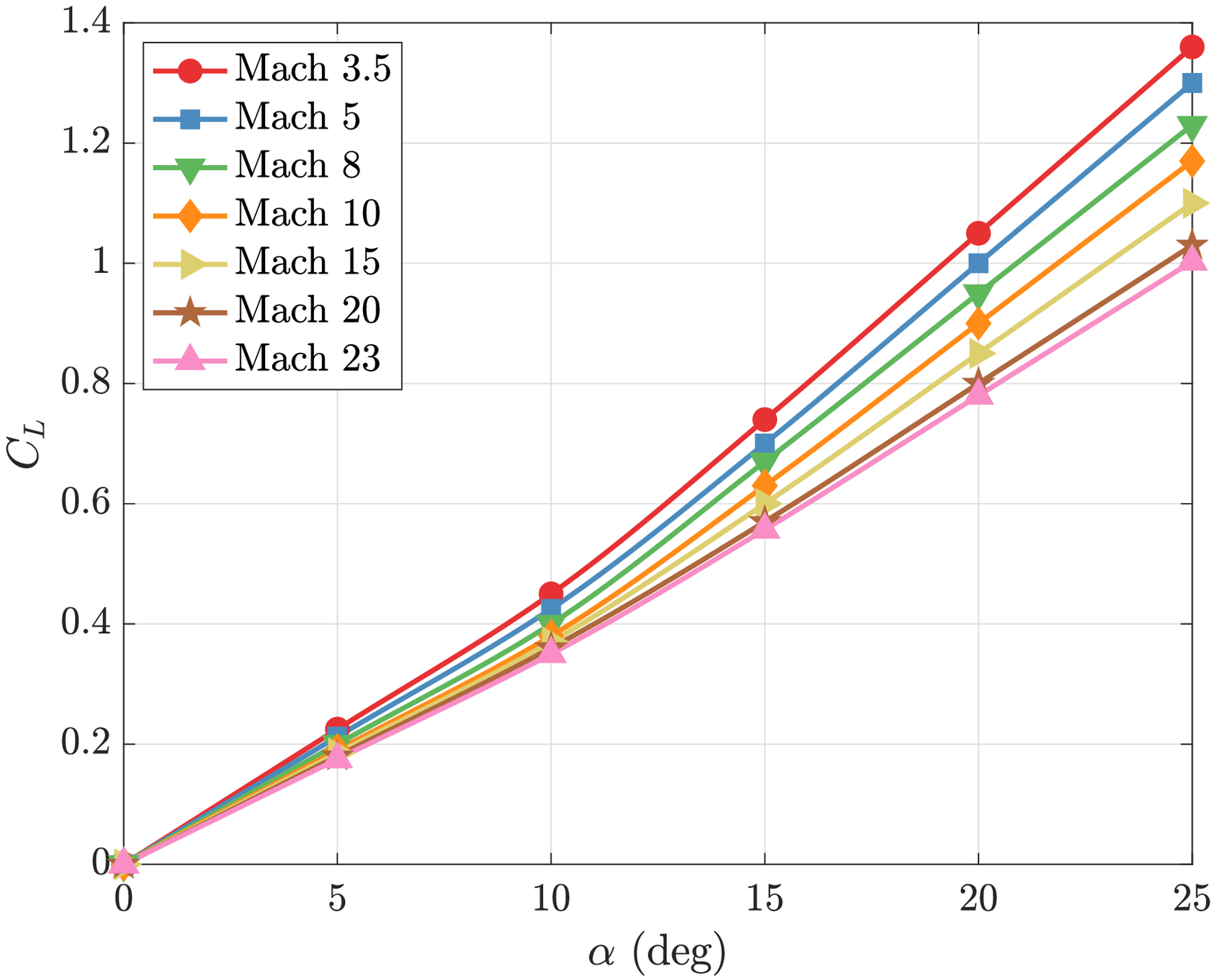}}
  &
  \subfloat[Drag coefficient, $C_D$ vs. angle of attack, $\alpha$.]{\includegraphics[width=.475\textwidth]{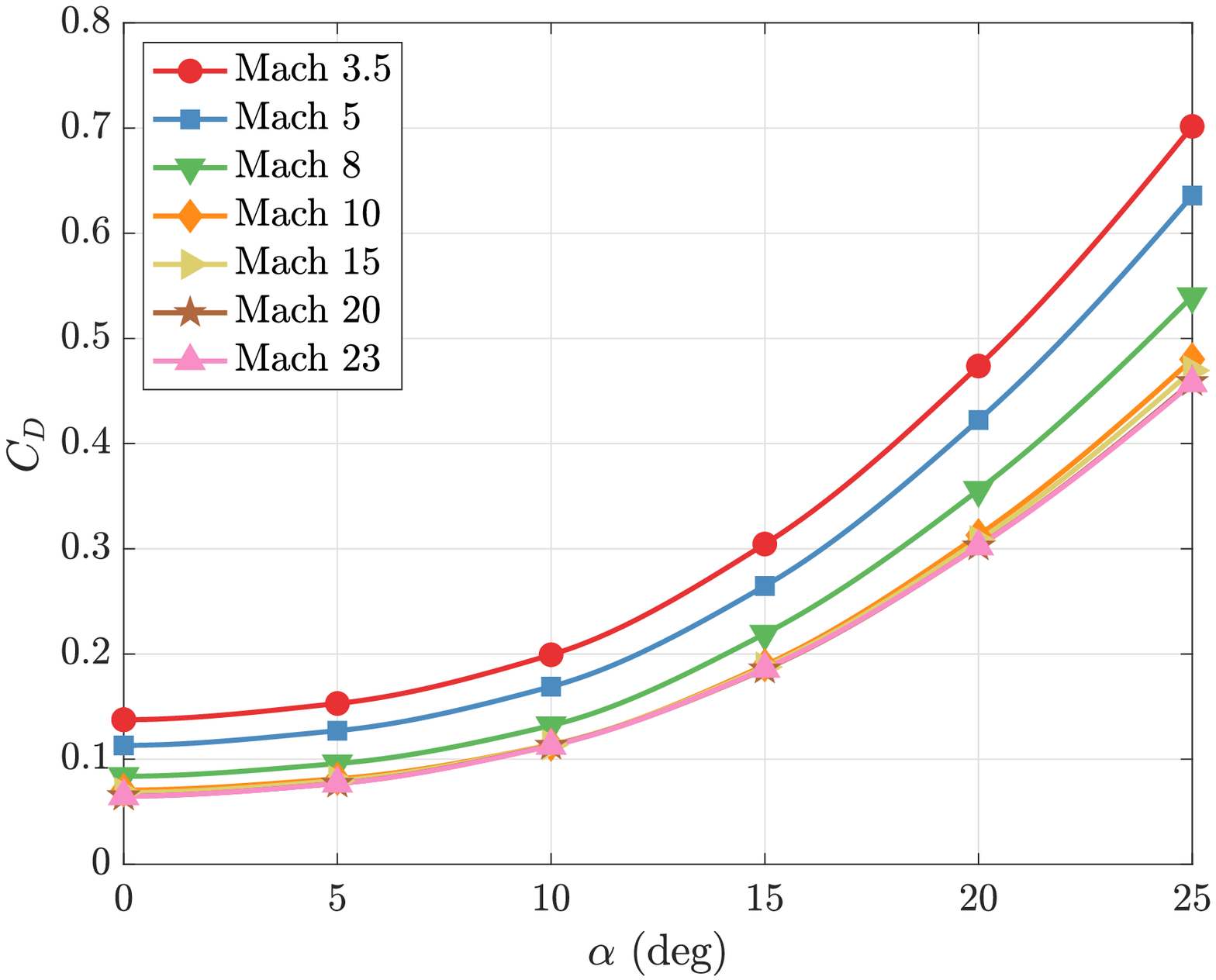}}
  \\
  \subfloat[Lift-to-drag ratio, $L/D$ vs. angle of attack, $\alpha$.]{\includegraphics[width=.475\textwidth]{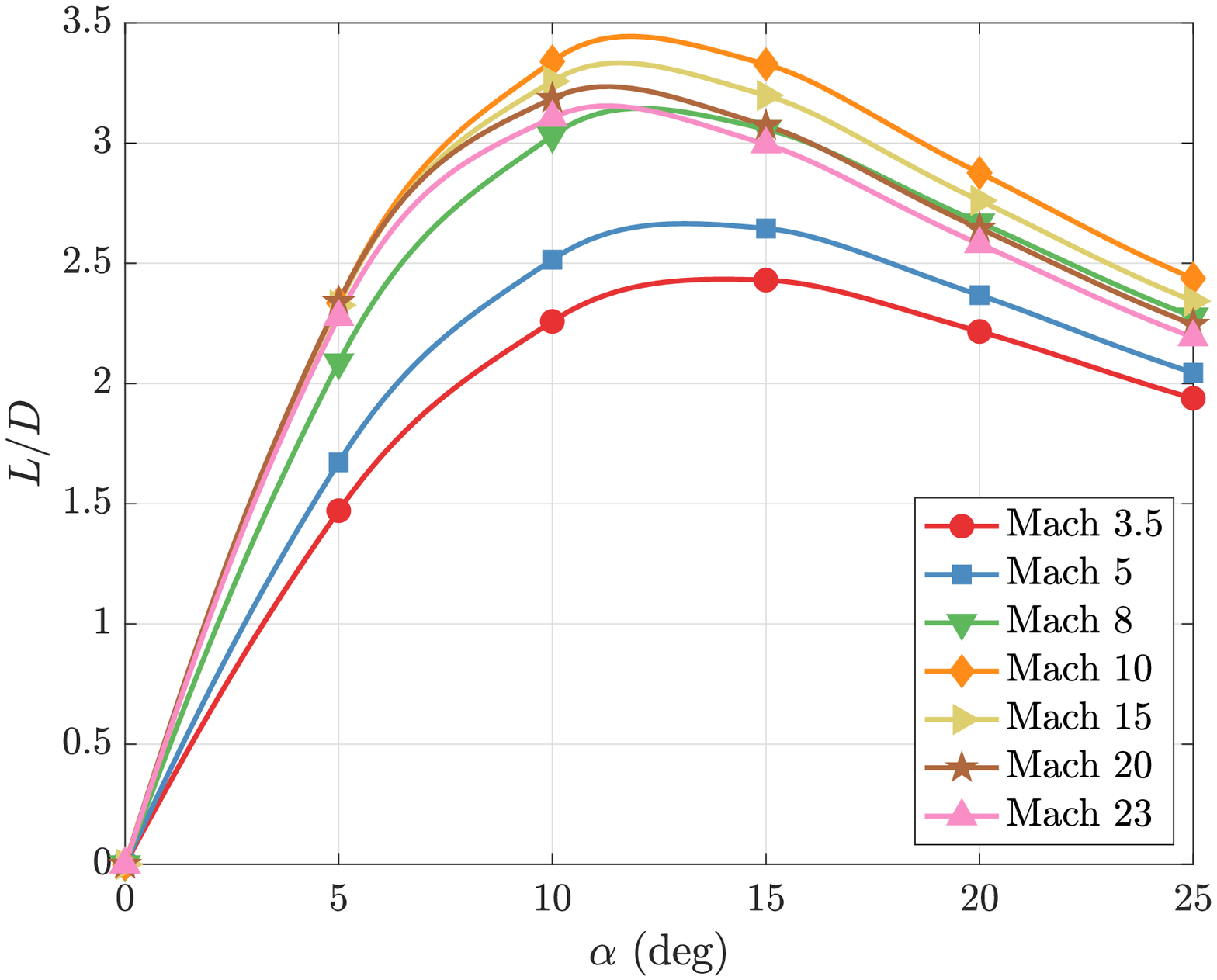}}
  &
  \subfloat[Drag polar model for Mach 3.5.]{\includegraphics[width=.475\textwidth]{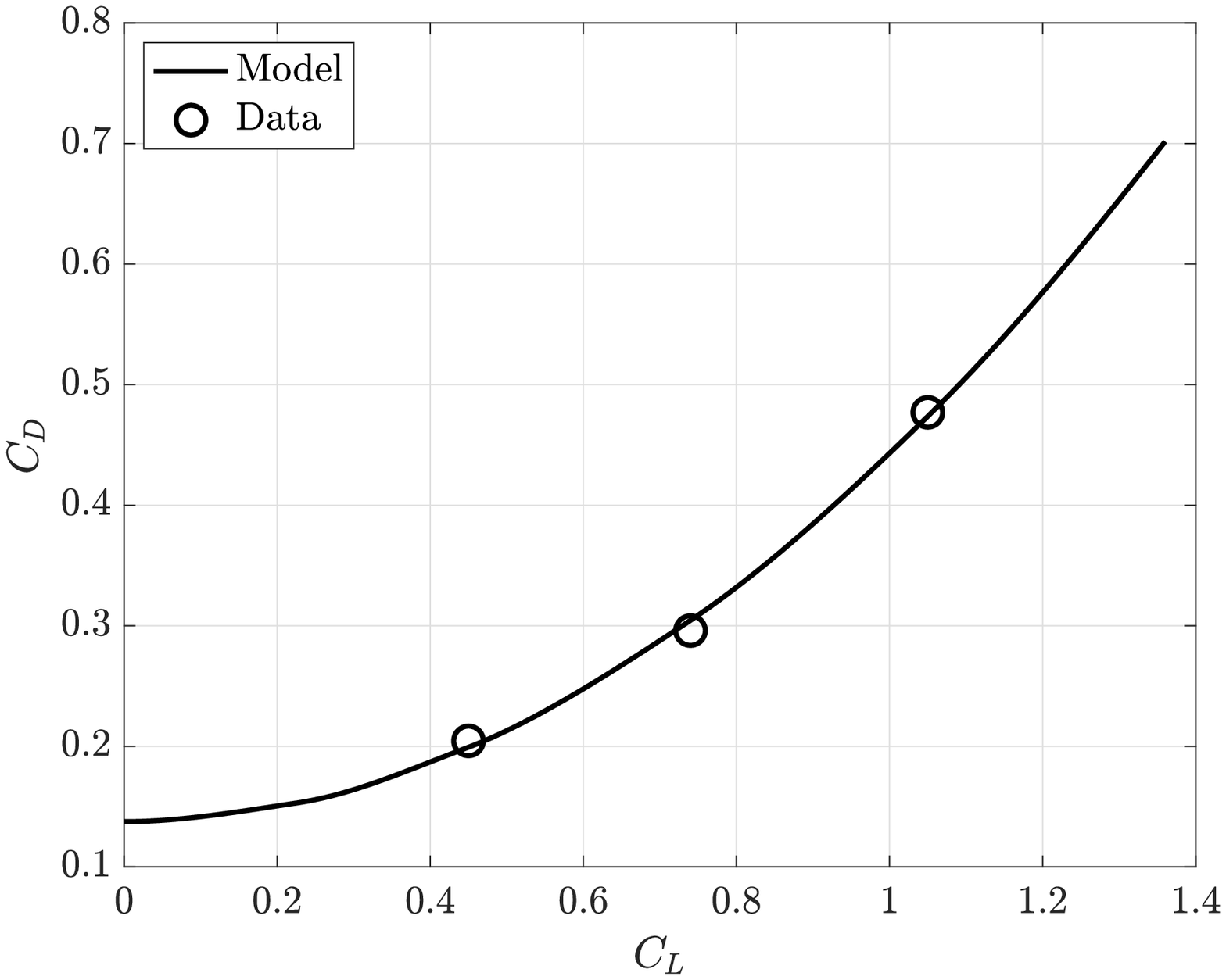}}
  \end{tabular}
  \caption{Aerodynamic model for the entry vehicle.}\label{fig:CAV-LD}
\end{figure}

The lift and drag coefficient data (at each Mach number for which data is available) is extended to include values at $\alpha = \{0, 5, 25\}~\textrm{deg}$ as follows.  First, it is assumed that $C_L = 0$ at zero angle of attack.  Linear interpolation of the $C_L$ data at $\alpha = \{0, 10\}$ is then employed to obtain the value of $C_L$ at $\alpha = 5$~deg.  Similarly, the $C_L$ data at $\alpha = \{15, 20\}$ is linearly extrapolated to obtain the value of $C_L$ at $\alpha = 25$~deg.  Next, when the Mach number is held constant, the value of $C_D$ is assumed to behave according to the drag polar model
\begin{equation}\label{eq:DragPolar}
\begin{array}{rcl}
C_D &=& C_{D0} + K C_L^2,
\end{array}
\end{equation}
where $C_{D0}$ is the zero-lift drag coefficient and $K$ is the drag polar parameter.  The values of the parameters $\{C_{D0}, K\}$ are obtained by employing a least squares fit to the $C_L$ and $C_D$ data available at the desired Mach number.  After attaining the optimal fit, Eq.~\eqref{eq:DragPolar} is then employed to obtain the values for $C_D$ at $\alpha = \{0, 5, 10, 15, 20, 25\}$~deg.  Finally, given the aforementioned data, 2D interpolation is employed to estimate values of the lift and drag coefficients as functions of angle of attack and Mach number\cite{makima}.  Figure~\ref{fig:CAV-LD} illustrates the aerodynamic model.

\begin{table}[h]
  \centering
  \caption{Entry Vehicle Data.\label{tab:CAV-data}}
  \renewcommand{\baselinestretch}{1}\normalsize\normalfont
 \begin{tabular}{lccl}\hline
    Description & Symbol & Units & Value \\\hline
    Mass & $m$ & kg & $907.186$ \\
    Reference Area & $S$ & $\textrm{m}^2$ & $0.48387$ \\\hline
  \end{tabular}
\end{table}

\clearpage 
\subsection{Equations of Motion}\label{sect:EOM}
\subsubsection{Non-Vertical Flight}\label{sect:EOM-geo}
Phases 2-7 do not contain any periods of vertical flight.  During these phases, the equations of motion are given by
\begin{equation}\label{eq:EOM-geo}
\begin{array}{rcl}
\dot{h} &=& v\sin\gamma, \\
\dot{\phi} &=& \frac{v}{r\cos\theta}\cos\gamma\sin\psi, \\
\dot{\theta} &=& \frac{v}{r}\cos\gamma\cos\psi, \\

\dot{v} &=& \frac{1}{m} \left(T\cos\alpha - D\right) -
                    \frac{\mu_e}{r^2}\sin\gamma + 
                    r\omega_e^2\cos\theta\left(
                    \sin\gamma\cos\theta - \cos\gamma\sin\theta\cos\psi
                    \right), \\

\dot{\gamma} &=& \frac{\cos\sigma}{mv}\left(T\sin\alpha + L\right) +
                                \cos\gamma\left(\frac{v}{r} - \frac{\mu_e}{r^2v}\right) + 
                                2\omega_e\cos\theta\sin\psi + 
                                \frac{r\omega_e^2}{v}\cos\theta\left(
                                \cos\gamma\cos\theta + \sin\gamma\sin\theta\cos\psi
                                \right), \\

\dot{\psi} &=& \frac{\sin\sigma}{mv\cos\gamma}\left(T\sin\alpha + L\right) +
                \frac{v}{r}\cos\gamma\sin\psi\tan\theta - 
                2\omega_e\left(\tan\gamma\cos\theta\cos\psi - \sin\theta\right) + 
                \frac{r\omega_e^2}{v\cos\gamma}\sin\theta\cos\theta\sin\psi, \\

\dot{m} & = & \frac{T}{I_{SP} g_0},
\end{array}
\end{equation}
where $h$ is the altitude of the vehicle above the spherical Earth, $\phi$ is the Earth-relative longitude, $\theta$ is the geocentric latitude, $v$ is the Earth-relative speed, $\gamma$ is the Earth-relative flight path angle, $\psi$ is the azimuth angle, $\alpha$ is the vehicle's angle of attack, and $\sigma$ is the bank angle.  It is noted that the bank angle defines the angle between the lift vector and the plane formed by the position and Earth-relative velocity vectors.  Next, the geocentric radius is $r = h + R_e$, where $R_e$ is the radius of the Earth, $\omega_e$ is the Earth's rotation rate, $\mu_e$ is the gravitational parameter of the Earth, and $g_0$ is the Earth's standard gravitational acceleration.  Finally, $m$ denotes the vehicle mass, $T$ and $I_{SP}$ are the thrust and specific impulse of the boost vehicle, and the relations
\begin{equation}\label{eq:LD}
\begin{array}{rcl}
L = q S C_L, \\
D = q S C_D,
\end{array}
\end{equation}
define the lift and drag force magnitudes, denoted by $L$ and $D$ respectively, where $S$ is the vehicle reference area, $C_L$ and $C_D$ are the lift and drag coefficients, $q = \rho v^2 / 2$ is the dynamic pressure, and $\rho$ is the ambient atmospheric density.

The equations of motion given by Eq.~\eqref{eq:EOM-geo} are written in a general form that is suitable for describing the motion in phases 2-7.  However, it is noted that some of the terms are equal to zero during particular phases of flight.  In particular, phases 4, 5, and 6 are considered exo-atmospheric phases.  Thus, the lift and drag force terms are zero.  Similarly, the thrust terms are zero in phases 5, 6, and 7 because these phases are non-propulsive.

Lastly, two additional differential equations are appended to Eq.~\eqref{eq:EOM-geo} in order to account for realistic rates of change in the angle of attack and bank angle.  The two equations are given by
\begin{equation}\label{eq:aoa-bank}
\begin{array}{rcl}
\dot{\alpha} &=& u_\alpha, \\
\dot{\sigma} &=& u_\sigma,
\end{array}
\end{equation}
where $u_\alpha$ and $u_\sigma$ are, respectively, the angle of attack rate and the bank angle rate.  Together, the variables $\{ h, \phi, \theta, v, \gamma, \psi, \alpha, \sigma, m \}$ comprise the state in phases 2-4 and $\{ h, \phi, \theta, v, \gamma, \psi, \alpha, \sigma\}$ is the state in phases 5-7.  Likewise, $\{ u_\alpha, u_\sigma \}$ is the control in phases 2-7.

\subsubsection{Vertical Flight}\label{sect:EOM-Euler}
Phases 1 and 8 begin and end, respectively, in vertical flight ($\gamma = \pm 90$ deg).  Notice that the azimuth rate of change in Eq.~\eqref{eq:EOM-geo} is undefined at the vertical flight condition.  Thus, the following alternative to Eq.~\eqref{eq:EOM-geo} is employed during phases 1 and 8 in order to remove the singularity during vertical flight.  The equations of motion are given by
\begin{equation}\label{eq:EOM-Euler}
\begin{array}{rcl}
\dot{h} &=& v\left(1 - 2\left(\e{2}^2 + \e{3}^2\right)\right), \\
\dot{\phi} &=& \frac{2v}{r\cos\theta}\left(\e{1}\e{2} + \e{3}\eta\right), \\
\dot{\theta} &=& \frac{2v}{r}\left(\e{1}\e{3} - \e{2}\eta\right), \\

\dot{v} &=& \frac{1}{m} \left(T\cos\alpha - D\right) -
                    \frac{\mu_e}{r^2}\left(1 - 2\left(\e{2}^2 + \e{3}^2\right)\right) + 
                    r\omega_e^2\cos\theta\left(
                    \cos\theta\left(1 - 2\left(\e{2}^2 + \e{3}^2\right)\right) - 
                    2\sin\theta\left(\e{1}\e{3} - \e{2}\eta\right)
                    \right),
\end{array}
\end{equation}
and
\begin{equation}\label{eq:EulerRates}
\begin{array}{rclcrcl}
\edot{1} &=& \frac{1}{2}\left(\phantom{-}\eta\w{1} - \e{3}\w{2} + \e{2}\w{3}\right),
& \quad & 
\edot{2} &=& \phantom{-}\frac{1}{2}\left(\e{3}\w{1} - \eta\w{2} - \e{1}\w{3}\right), \\
\edot{3} &=& \frac{1}{2}\left(-\e{2}\w{1} + \e{1}\w{2} + \eta\w{3}\right),
& \quad &
\dot{\eta} &=& -\frac{1}{2}\left(\e{1}\w{1} + \e{2}\w{2} + \e{3}\w{3}\right),
\end{array}
\end{equation}
where $\dot{m} = T/\left(I_{SP} g_0\right)$ and $\dot{\alpha} = u_\alpha$ still hold, and where
\begin{equation}\label{eq:w2w3}
\begin{array}{rcl}
\w{2} &=& -\frac{1}{mv}T\sin\alpha - 
                 2\left(\frac{v}{r} - \frac{\mu_e}{r^2v}\right)
                 \left(\e{1}\e{3} + \e{2}\eta\right) - 
                 4\w{e}\left[
                 \sin\theta\left(\e{1}\e{2} - \e{3}\eta\right) + 
                 \cos\theta\left(\e{2}\e{3} + \e{1}\eta\right)
                 \right] \\
                 & & - 
                 \frac{2r\w{e}^2}{v}\cos\theta\left[
                 \cos\theta\left(\e{1}\e{3} + \e{2}\eta\right) - 
                 \sin\theta\left(\frac{1}{2} - \e{1}^2 - \e{2}^2\right)
                 \right], \\
                 
\w{3} &=& \frac{1}{mv}\left(T\sin\alpha + L\right) + 
                  2\left(\frac{v}{r} - \frac{\mu_e}{r^2v}\right)
                  \left(\e{1}\e{2} - \e{3}\eta\right) - 
                  4\w{e}\left[
                  \sin\theta\left(\e{1}\e{3} + \e{2}\eta\right) + 
                  \cos\theta\left(\frac{1}{2} - \e{1}^2 - \e{2}^2\right)
                  \right] \\
                  & & + 
                  \frac{2r\w{e}^2}{v}\cos\theta\left[
                  \cos\theta\left(\e{1}\e{2} - \e{3}\eta\right) - 
                  \sin\theta\left(\e{2}\e{3} + \e{1}\eta\right)
                  \right].
\end{array}
\end{equation}
It is noted that the variables $\{\e{1}, \e{2}, \e{3}, \eta\}$ are Euler parameters (unit quaternion), $\{\w{1}, \w{2}, \w{3}\}$ are angular velocity components ($w_1$ is treated as a control component), and all other variables retain their definitions from Section~\ref{sect:EOM-geo}.  For further information regarding Eqs.~\eqref{eq:EOM-Euler}--\eqref{eq:w2w3} the reader is referred to Ref.~\cite{Miller2}.  Finally, it is noted that the variables $\{ h, \phi, \theta, v, \e{1}, \e{2}, \e{3}, \eta, \alpha, m \}$ comprise the state in phase 1, $\{ h, \phi, \theta, v, \e{1}, \e{2}, \e{3}, \eta, \alpha \}$ is the state in phase 8, and $\{ u_\alpha, w_1 \}$ is the control in phases 1 and 8.

\subsubsection{State and Control Conversions}\label{sect:convert-state-and-control}
The following mappings will prove useful throughout the remainder of this paper.  First, the Euler parameters $\{\e{1}, \e{2}, \e{3}, \eta\}$ of Section~\ref{sect:EOM-Euler} define the flight path angle $\gamma$, azimuth angle $\psi$, and bank angle $\sigma$ as
\begin{equation}\label{eq:convert-state}
\begin{array}{rcl}
\gamma &=& \tan^{-1}\left( 
                       \frac{1}{2} - \e{2}^2 - \e{3}^2~~,~~
                       \sqrt{(\e{1}^2 + \eta^2)(\e{2}^2 + \e{3}^2)}
                       \right), \\
\psi &=& \tan^{-1}\left(
               \e{1}\e{2} + \e{3}\eta~~,~~
               \e{1}\e{3} - \e{2}\eta
               \right), \\
\sigma &=& \tan^{-1}\left(
                    -\e{3}\e{1} - \e{2}\eta~~,~~
                    \e{2}\e{1} - \e{3}\eta
                    \right),
\end{array}
\end{equation}
respectively, where $\tan^{-1}(\cdot,\cdot)$ is the four-quadrant inverse tangent operator.  Next, the angular velocity component $\w{1}$ is related to the bank angle rate $u_{\sigma}$ by
\begin{equation}\label{eq:convert-control}
u_{\sigma} = \w{1} - 
\frac{\frac{1}{2} - \e{2}^2 - \e{3}^2}{(\e{1}^2 + \eta^2)(\e{2}^2 + \e{3}^2)}
\left[
\w{2}(\e{2}\e{1} - \e{3}\eta) + \w{3}(\e{3}\e{1} + \e{2}\eta)
\right],
\end{equation}
where $\w{2}$ and $\w{3}$ are defined in Eq.~\eqref{eq:EOM-Euler}.  Observing Eq.~\eqref{eq:convert-control}, it is noticed that $u_{\sigma} = \w{1}$ during horizontal flight ($\gamma = 0~\textrm{deg}$ or equivalently $\e{2}^2 + \e{3}^2 = \frac{1}{2}$) or when the Earth-relative velocity direction is stationary as viewed by an observer in the LVLH frame ($\w{2} = \w{3} = 0$).

\section{Problem Formulation \label{sect:ProblemFormulation}}
The combined ascent-entry trajectory optimization problem is stated as an optimal control problem as follows.  Determine the state and control signals as well as the unknown initial and terminal times in each phase such that the performance index of Section~\ref{sect:Cost} is minimized.  The minimizing solution must simultaneously satisfy the vehicle dynamics of Section~\ref{sect:EOM} as well as the path constraints, heat load constraint, boundary conditions, and interior point constraints of Sections~\ref{sect:Con-Path}--\ref{sect:Con-Interior}.  Detailed descriptions of the objective and constraints involved in the problem are provided next.

\subsection{Performance Index}\label{sect:Cost}
A primary goal of this research is to produce ascent-entry trajectories that are suitable reference trajectories for outer-loop guidance applications.  The reference trajectory could be generated once for the entire mission, or it could be updated periodically as flight progresses.  In either case, it is desirable to produce reference trajectories with wide control margins.  That way perturbations during actual flight can be counteracted using the remaining control authority.  With the aforementioned goal in mind, the performance index employed in this research is designed to produce trajectories with wide control margins in all phases of flight.  The performance index is given by
\begin{equation}\label{eq:Cost-total}
\C{J} = \sum_{p = 1}^{8} \C{J}^{(p)},
\end{equation}
where $\C{J}^{(p)}$ is the integrated cost across phase $p$.  The general form of $\C{J}^{(p)}$ in each phase is given by
\begin{equation}\label{eq:Cost-phase}
\C{J}^{(p)} = \int_{t_0^{(p)}}^{t_f^{(p)}} \left[
\left( \frac{\alpha - \bar{\alpha}}{\alpha_{\max}} \right)^2 + 
\left( \frac{u_\alpha}{u_{\alpha,\max}} \right)^2 + 
\left( \frac{u_\sigma}{u_{\sigma,\max}} \right)^2 \right] dt
\end{equation}
where $t_0^{(p)}$ and $t_f^{(p)}$ are the initial and terminal times of phase $p$ and $\{ \bar{\alpha}, \alpha_{\max}, u_{\alpha,\max}, u_{\sigma,\max} \}$ are phase-dependent constants. 

Equation~\eqref{eq:Cost-phase} is employed in all eight phases of flight with the following exceptions.  First, during the exo-atmospheric coasting phases the first term of the integrand in Eq.~\eqref{eq:Cost-phase} is removed.  Second, the control component $\w{1}$ replaces $u_{\sigma}$ in Eq.~\eqref{eq:Cost-phase} during phases 1 and 8.  Finally, Section~\ref{sect:sigStudy} studies the effects of adding an additional penalty term to the integrand of Eq.~\eqref{eq:Cost-phase} during phase 8, noting that the penalty term is intended to reduce phugoid oscillations during entry.  The precise definition of the additional penalty term and its effects on generated trajectories is covered later in Section~\ref{sect:sigStudy}.

\subsection{Path Constraints}\label{sect:Con-Path}
\subsubsection{Bounds on State and Control Variables}
Limits on the state and control variables in each phase are summarized as follows.  First, the angle of attack is bounded during ascent (phases 1-4) by
\begin{equation}\label{eq:bounds-aoa-boost}
|\alpha| \leq \alpha_{\max},
\end{equation}
and during entry (phases 7 and 8) by
\begin{equation}\label{eq:bounds-aoa-entry}
0 \leq \alpha \leq \alpha_{\max},
\end{equation}
noting that the angle of attack is free during the exo-atmospheric coasting phases (phases 5 and 6).  Next, the angle of attack and bank angle rates are bounded in all phases by
\begin{equation}\label{eq:bounds-uaub}
\begin{array}{rcl}
|u_\alpha| & \leq & u_{\alpha,\max}, \\
|u_\sigma| & \leq & u_{\sigma,\max},
\end{array}
\end{equation}
where again it is noted that $w_1$ replaces $u_\sigma$ in phases 1 and 8.  Finally, phases 4-6 are exo-atmospheric and must satisfy the minimum altitude constraint
\begin{equation}\label{eq:bounds-altitude-exo}
h \geq h_{\textrm{atm}},
\end{equation} 
where $h_{\textrm{atm}}$ is a reasonably chosen altitude above which the Earth's atmospheric effects are considered small.

\subsubsection{Sensed Acceleration}
During entry (phases 7 and 8), the glide vehicle must maintain a sufficiently low sensed acceleration.  The sensed acceleration (expressed in g's) is defined as
\begin{equation}\label{eq:g-force}
n = \frac{1}{m g_0}\sqrt{L^2 + D^2},
\end{equation}
where $g_0$ is the standard gravitational acceleration of the Earth given in Table~\ref{tab:Earth-data}.  The path constraint is then applied in phases 7 and 8 as
\begin{equation}\label{eq:maxG}
n \leq n_{\max},
\end{equation}
where $n_{\max}$ is the specified upper limit on sensed acceleration.

\subsubsection{Dynamic Pressure}
Maximum dynamic pressure during ascent occurs in phase 1.  As such, the path constraint
\begin{equation}\label{eq:maxQ}
q \leq q_{\max},
\end{equation} 
is enforced in phase 1 to ensure structural integrity of the boost vehicle.  Two additional dynamic pressure constraints are applied in phases 7 and 8 in order to separate entry flight into low and high dynamic pressure phases.  The dynamic pressure constraints are given by
\begin{equation}\label{eq:minQ-1}
q \leq q_{\min},
\end{equation} 
in phase 7 and
\begin{equation}\label{eq:minQ-2}
q \geq q_{\min},
\end{equation}
in phase 8, where $q_{\min}$ specifies the boundary between the low and high dynamic pressure phases.  Together, Eqs.~\eqref{eq:minQ-1} and \eqref{eq:minQ-2} force the entry vehicle to fall into the Earth's atmosphere without skipping back out and losing aerodynamic control.

\subsubsection{Heating Rate}
Thermal protection of the entry vehicle requires that a sufficiently low stagnation point heating rate be maintained during entry (phases 7 and 8).  The stagnation point heating rate is computed via the Chapman equation\cite{Detra1} as
\begin{equation}\label{eq:heatRate}
\dot{Q} = \kappa
                \left( \frac{\rho}{\rho0} \right)^{0.5}
                \left( \frac{v}{v_c} \right)^{3.15},
\end{equation}
where $\kappa = 199.87~\textrm{MW}/\textrm{m}^2$, $\rho_0 = 1.225~\textrm{kg}/\textrm{m}^3$, and $v_c = \sqrt{\mu_e/R_e} = 7.9053~\textrm{km}/\textrm{s}$.  The heating rate path constraint is then given as
\begin{equation}\label{eq:heatRateLimit}
\dot{Q} \leq \dot{Q}_{\max},
\end{equation}
where $\dot{Q}_{\max}$ is the specified upper limit.

\subsection{Heat Load \label{sect:Con-Integral}}

In addition to limits placed on the heating rate during entry, thermal protection of the entry vehicle also requires that the heating load be sufficiently small during entry.  The heating load experienced during entry is defined by the integral
\begin{equation}\label{eq:Q}
Q = \int_{t_0^{(7)}}^{t_f^{(8)}} \dot{Q} dt,
\end{equation}
where $t_0^{(7)}$ and $t_f^{(8)}$ are the initial and terminal times, respectively, of phases 7 and 8, and where $\dot{Q}$ is given in Eq.~\ref{eq:heatRate}.  Thus, the constraint
\begin{equation}\label{eq:Qmax}
  Q \leq Q_{\max}
\end{equation}
is imposed during atmospheric entry.  

\subsection{Initial and Terminal Conditions}\label{Con-BCs}

The initial and terminal conditions are listed in Table~\ref{tab:BCs}.  The initial conditions are obtained by propagating a 1D simulation of the rocket (with $\alpha = 0$) from launch ($t = 0$) until the tower has been safely cleared (taken to be $50~\textrm{m}$ above the launch pad).  It is noted that the initial conditions correspond to vertically upwards flight with the belly of the rocket facing West.  Next, the terminal boundary conditions specify the requirements for target impact.  The impact requirements employed here correspond to vertically downwards flight and an impact speed of $1.219~\textrm{km/s}$.

\begin{table}[h]
  \centering
  \caption{Initial and Terminal Conditions.\label{tab:BCs}}
  \renewcommand{\baselinestretch}{1}\normalsize\normalfont
 \begin{tabular}{lccrr}\hline
     Description & Symbol & Units & Initial & Terminal \\\hline
    Time & $t$ & s & 2.52 & FREE \\
    Altitude & $h$ & km &  0.167 & 0.000 \\
    Longitude & $\phi$ & deg & -120.63 & -192.30 \\
    Latitude & $\theta$ & deg & 34.58 & 8.70 \\
    Speed & $v$ & km/s & 0.040 & 1.219 \\
    Euler Parameter & $\e{1}$ & - & 0 & 0 \\
    Euler Parameter & $\e{2}$ & - & 0 & FREE \\
    Euler Parameter & $\e{3}$ & - & 0 & FREE \\
    Euler Parameter & $\eta$ & - & 1 & 0 \\
    Angle of Attack & $\alpha$ & deg & 0 & 0 \\
    Mass & $m$ & kg & 85743 & - \\\hline
  \end{tabular}
\end{table}

\subsection{Interior Point Constraints}\label{sect:Con-Interior}
All eight phases occur sequentially with both state and time continuity maintained from one phase to the next.  Thus, the following interior point constraints are employed to maintain continuity at each phase boundary.  First, let the time domain of phase $p$ be denoted $[t_0^{(p)},t_f^{(p)}]$, where $t_0^{(p)}$ and $t_f^{(p)}$ are the initial and terminal times of phase $p$.  It follows that the phase transitions occur at $t_0^{(p)},~p = 2,\ldots,8$ or equivalently at $t_f^{(p)},~p = 1,\ldots,7$.  Given that the time sequence of each boost phase is known, the time continuity constraints
\begin{equation}
\begin{array}{rcccl}
t_f^{(1)} &=& t_0^{(2)} &=& t_{S1}, \\
t_f^{(2)} &=& t_0^{(3)} &=& t_{S2}, \\
t_f^{(3)} &=& t_0^{(4)} &=& t_{\textrm{fairing}}, \\
t_f^{(4)} &=& t_0^{(5)} &=& t_{S3},
\end{array}
\end{equation} 
enforce the appropriate phase transition times, noting that $t_{\textrm{fairing}}$ is the time at which fairing separation occurs and that $\{ t_{S1}, t_{S2}, t_{S3} \}$ are the times at which engine burnout occurs, respectively, for stages 1-3.  Each of the remaining phase transition times are free variables in the problem.  Thus, time continuity is enforced by requiring
\begin{equation}\label{eq:continuity-time}
t_0^{(p+1)} - t_f^{(p)} = 0, \quad p = 5,\ldots,7. \\
\end{equation}

Similar to the time continuity constraints, state continuity is maintained at each phase boundary via the interior point constraints
\begin{equation}\label{eq:continuity-state}
\m{y}\left({t_0^{(p+1)}}\right) - \m{y}\left({t_f^{(p)}}\right) = \m{0},~p = 1,\ldots,7,
\end{equation}
where $\m{y} = [h,\phi,\theta,v,\gamma,\psi,\alpha,\sigma]\tr$, noting that Eq.~\eqref{eq:convert-state} is employed to convert the Euler parameters in phases 1 and 8 to $\gamma$, $\psi$, and $\sigma$ at $t_f^{(1)}$ and $t_0^{(8)}$.  In addition, the unit norm constraint
\begin{equation}
\left( \e{1}^2 + \e{2}^2 + \e{3}^2 + \eta^2 \right)\big\rvert_{t_0^{(8)}} = 1,
\end{equation}
ensures the Euler parameters have unit norm in phase 8.  The unit norm constraint is necessary because the values of the Euler parameters are not fully specified at either endpoint of phase 8.  

Next, it is noted that the state continuity constraints of Eq.~\ref{eq:continuity-state} do not include mass.  Instead, the discontinuous change in mass of the boost vehicle at stage and fairing separation is handled by employing the mass transition constraints
\begin{equation}
\begin{array}{rcl}
m\left({t_0^{(2)}}\right) &=& m_{S2}, \\
m\left({t_0^{(3)}}\right) &=& m_{S3}, \\
m\left({t_0^{(4)}}\right) - m\left({t_f^{(3)}}\right) &=& m_{\textrm{fairing}},
\end{array}
\end{equation} 
where $m_{\textrm{fairing}}$ is the mass of the fairing and where $m_{S2}$ and $m_{S3}$ denote the total mass of the boost vehicle at stage 2 and stage 3 ignition, respectively.  Finally, it is noted that mass is only a component of the state in phases 1-4.  Thus, no further mass transition constraints are necessary for the exo-atmospheric coast and entry phases.

Next, the following interior point constraints define particular phase boundary conditions that must be satisfied.  First, $t_f^{(5)} = t_0^{(6)}$ is the point at which the peak altitude is reached and where the entry vehicle separates from the final stage of the boost vehicle.  Thus, the boundary between phases 5 and 6 is defined by the payload separation conditions
\begin{equation}
\begin{array}{rcccl}
\alpha\left({t_f^{(5)}}\right) &=& \alpha\left({t_0^{(6)}}\right) &=& 0, \\
\sigma\left({t_f^{(5)}}\right) &=& \sigma\left({t_0^{(6)}}\right) &=& 0,
\end{array}
\end{equation}
and the peak altitude conditions
\begin{equation}
\begin{array}{rcccl}
h\left({t_f^{(5)}}\right) &=& h\left({t_0^{(6)}}\right) & \in & [h_{\textrm{peak},\min}~,~h_{\textrm{peak},\max}], \\
\gamma\left({t_f^{(5)}}\right) &=& \gamma\left({t_0^{(6)}}\right) &=& 0,
\end{array}
\end{equation} 
where $[h_{\textrm{peak},\min}~,~h_{\textrm{peak},\max}]$ defines the window of allowable peak altitudes.  Finally, the pierce point condition
\begin{equation}
h\left({t_f^{(6)}}\right) = h\left({t_0^{(7)}}\right) = h_{\textrm{atm}},
\end{equation}
defines the boundary between phase 6 (exo-atmospheric) and phase 7 (endo-atmospheric).

\subsection{Summary of Constraint Constants}\label{sect:CVV}
Tables~\ref{tab:CVV-1} and \ref{tab:CVV-2} summarize the values of the constraint constants defined throughout Section~\ref{sect:ProblemFormulation}.  It is noted that both $\dot{Q}_{\max}$ and $Q_{\max}$ are varied in the performance studies of Section~\ref{sect:results}, so their values are not provided here.  It is also noted that the value for $q_{\max}$ is calculated using Eq.~\ref{eq:maxQ} with the altitude and speed at maximum dynamic pressure given by Ref.~\cite{Min456UsersGuide}.

\begin{table}[h]
  \centering
  \caption{Values of Phase-Dependent Constraint Constants.\label{tab:CVV-1}}
  \renewcommand{\baselinestretch}{1}\normalsize\normalfont
 \begin{tabular}{lcccc}\hline
     Symbol & Units & Phases 1-4 & Phases 5-6 & Phases 7-8 \\\hline
    $\bar{\alpha}$ & deg & 0 & - & 11.86 \\
    $\alpha_{\max}$ & deg & 25 & - & 25 \\
    $u_{\alpha,\max}$ & deg/s & 10 & 10 & 10 \\
    $u_{\sigma,\max}$ & deg/s & 30 & 30 & 30 \\\hline
  \end{tabular}
\end{table}

\begin{table}[h]
  \centering
  \caption{Values of Phase-Independent Constraint Constants.\label{tab:CVV-2}}
  \renewcommand{\baselinestretch}{1}\normalsize\normalfont
 \begin{tabular}{lcc}\hline
     Symbol & Units & Value \\\hline
    $h_{\textrm{atm}}$ & km & 80 \\
    $h_{\textrm{peak},\min}$ & km & 100 \\
    $h_{\textrm{peak},\max}$ & km & 200 \\
    $m_{\textrm{fairing}}$ & kg & 400 \\
    $m_{S2}$ & kg & 38780 \\
    $m_{S3}$ & kg & 11110 \\
    $n_{\textrm{max}}$ & g & 12 \\
    $q_{\max}$ & kPa & 126.3 \\
    $q_{\min}$ & kPa & 12 \\
    $t_{\textrm{fairing}}$ & s & 179.1 \\
    $t_{S1}$ & s & 56.4 \\
    $t_{S2}$ & s & 117.1 \\
    $t_{S3}$ & s & 189.1 \\\hline
  \end{tabular}
\end{table}

\newpage 
\section{Results\label{sect:results}}
This section is divided into two sequential parts.  First, in Section~\ref{sect:sigStudy} a nominal solution to the ascent-entry trajectory optimization problem described in Section~\ref{sect:ProblemFormulation} is obtained where both the maximum stagnation point heating rate and heating load are unconstrained ($\dot{Q}_{\max} = Q_{\max} = \infty$).  The behavior of the nominal solution is then studied as an additional penalty term is applied to reduce phugoid oscillations in the generated trajectory.  Next, the results of Section~\ref{sect:sigStudy} are employed as a baseline of comparison for the studies carried out in Section~\ref{sect:conStudy}.  In particular, the nominal solution is compared against solutions obtained under varying maximum heating rate and heating load requirements.

Throughout Sections~\ref{sect:sigStudy} and \ref{sect:conStudy}, the following software and hardware are used.  All results are obtained using the the MATLAB optimal control software $\mathbb{GPOPS-II}$ \cite{Patterson2014}, noting that $\mathbb{GPOPS-II}$ employs $hp$-adaptive Legendre Gauss Radau collocation.  Next, the NLP solver employed is IPOPT \cite{Biegler2} set in full Newton mode at an accuracy tolerance of $10^{-6}$, and with both first and second derivatives supplied by the sparse central differencing method of Ref.~\cite{Patterson2012}.  The NLP max iteration count is set to $500$.  In addition, mesh refinement is carried out using the method of Ref.~\cite{Darby3} with a mesh error tolerance of $10^{-4}$ and the maximum number of mesh refinement iterations set to $10$.  Finally, all computations were performed on a 2.4 GHz 8-Core Intel Core i9 MacBook Pro running macOS Catalina version 10.15.7 with 32 GB of 2400 MHz DDR4 RAM and MATLAB version R2019b.

Due to the fact that a collocation method is employed and noting that Eq.~\ref{eq:EulerRates} consists of four differential equations with two degrees of freedom, special precautions must be taken to ensure that the system of equations generated in the numerical approximation of Eq.~\eqref{eq:EulerRates} is consistent.  Perfect integration of Eq.~\eqref{eq:EulerRates} would not result in any inconsistencies.  However, it is expected that the numerical approximation of Eq.~\eqref{eq:EulerRates} will have inconsistencies roughly on the same order of magnitude as the NLP solver tolerance.  Thus, a simple solution is to update Eq.~\eqref{eq:EulerRates} such that
\begin{equation}
\begin{array}{rcl}
\edot{2} &=& \e{3}\w{1} + \eta\w{2} - \e{1}\w{3} + u_1, \\
\edot{3} &=& -\e{2}\w{1} + \e{1}\w{2} + \eta\w{3} + u_2,
\end{array}
\end{equation}
where $\{u_1,u_2\}$ are slack variables with upper and lower bounds chosen to be one order of magnitude larger than the NLP solver tolerance.

\subsection{Reducing Phugoid Oscillations}\label{sect:sigStudy}
The nominal solution behavior is studied here where both the heating rate and heating load limits are relaxed to infinity (unconstrained).  As will become clear shortly, trajectories obtained using the performance index of Eq.~\eqref{eq:Cost-phase} tend to include phugoid oscillations during entry.  Phugoid oscillations are undesirable from a trajectory shaping perspective because each descent into the Earth's atmosphere typically coincides with large spikes in the sensed acceleration, dynamic pressure, and heating rate.  Instead, trajectories are preferred which remove the phugoid oscillations, providing more glide-like behavior and avoiding sudden spikes in the sensed acceleration, dynamic pressure, and heating rate.  

\begin{figure}[hbt!]
  \centering
  \includegraphics[width=.475\textwidth]{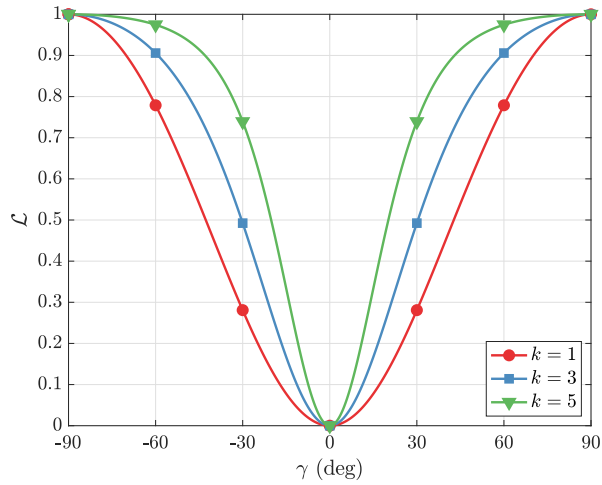}
  \caption{Phugoid oscillation penalty term as a function of flight path angle.}\label{fig:sigK}
\end{figure}

The following study explores the nominal solution behavior for the unconstrained heating rate and heating load case ($\dot{Q}_{\max} = Q_{\max} = \infty$) with and without the addition of a penalty term added to the integrand of Eq.~\ref{eq:Cost-phase} in phase 8.  The penalty term, denoted $\C{L}$, is defined as
\begin{equation}\label{eq:sig}
\C{L} = \left(\frac{C}{1 + e^{-k \sin\gamma}} - \frac{C}{2}\right)^2,
\end{equation}
where $\sin\gamma$ is expressed by the Euler parameters as $\sin\gamma = 1 - 2 \left( \e{2}^2 + \e{3}^2 \right)$, $k > 0$ is a design variable, and
\begin{equation}\label{eq:Cdef}
C = 2 \frac{1 + e^{-k}}{1 - e^{-k}},
\end{equation}
is a constant chosen such that $\C{L} = 1$ at $\gamma = \pm 90~\textrm{deg}$.  Figure~\ref{fig:sigK} illustrates the penalty term as a function of the flight path angle for the values of $k$ under study ($k = \{1, 3, 5\}$), noting that $k < 1$ produces nearly identical penalty profiles to $k = 1$ and that convergence issues in the resulting NLP begin to arise for values of $k \gg 5$.  Clearly, the penalty term is designed to reduce phugoid oscillations by incentivizing entry trajectories with a flight path angle near zero for the bulk of entry.  Thus, the penalty is zero at $\gamma = 0~\textrm{deg}$ and the penalty monotonically increases to one as the flight path angle increases or decreases from $\gamma = 0~\textrm{deg}$ to $\gamma = \pm 90~\textrm{deg}$.  The unit magnitude of the penalty term at $\gamma = \pm 90~\textrm{deg}$ is proportionate to the other terms in the integrand of Eq.~\eqref{eq:Cost-phase} so as to strike a balance between reducing phugoid oscillations and maintaining wide control margins in the generated trajectories.

\begin{figure}[hbt!]
  \centering
  \begin{tabular}{lr}
  \subfloat[Flight path angle, $\gamma(t)$ vs. time, $t$.]{\includegraphics[width=.475\textwidth]{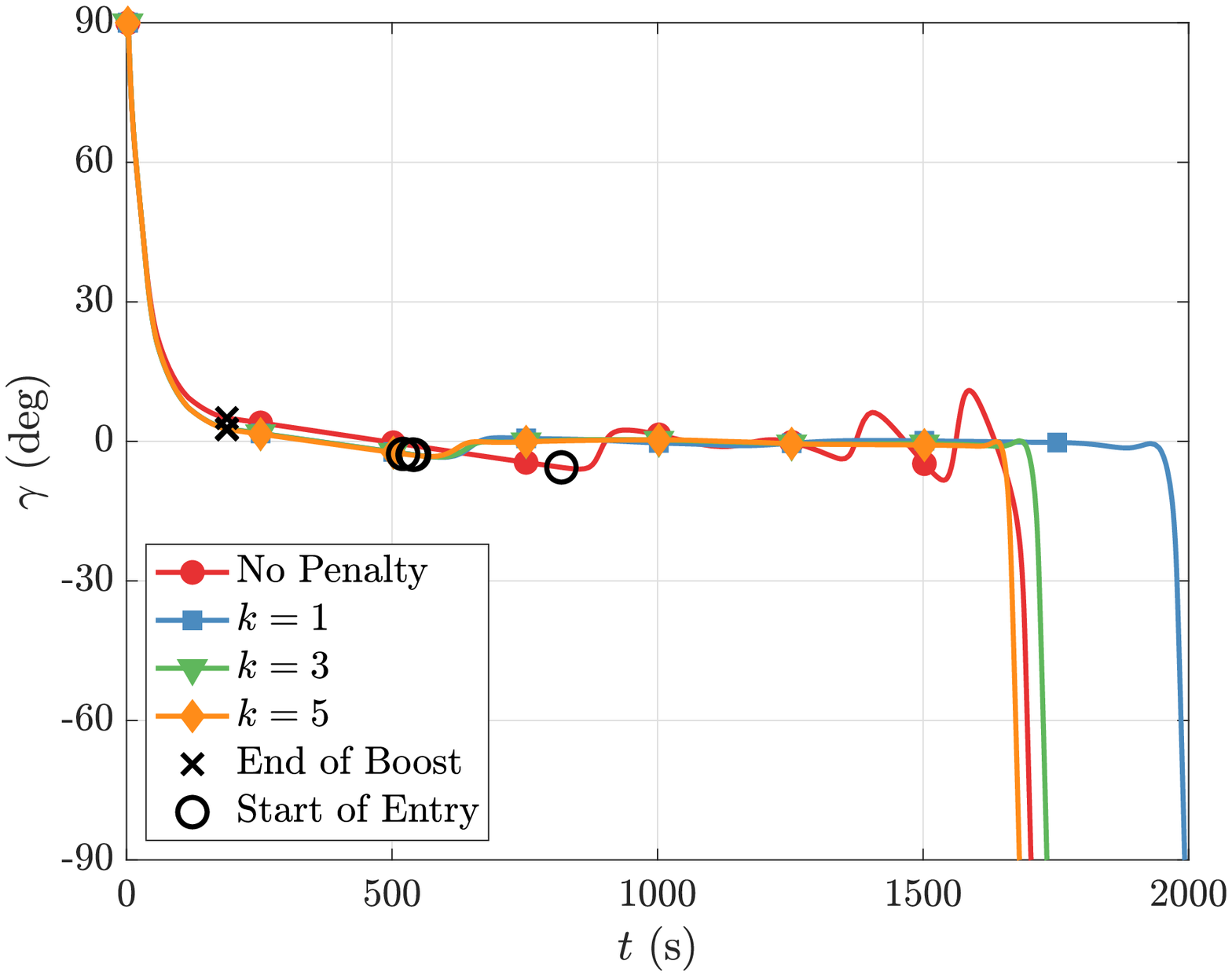}}
  &
  \subfloat[Enlarged view.]{\includegraphics[width=.475\textwidth]{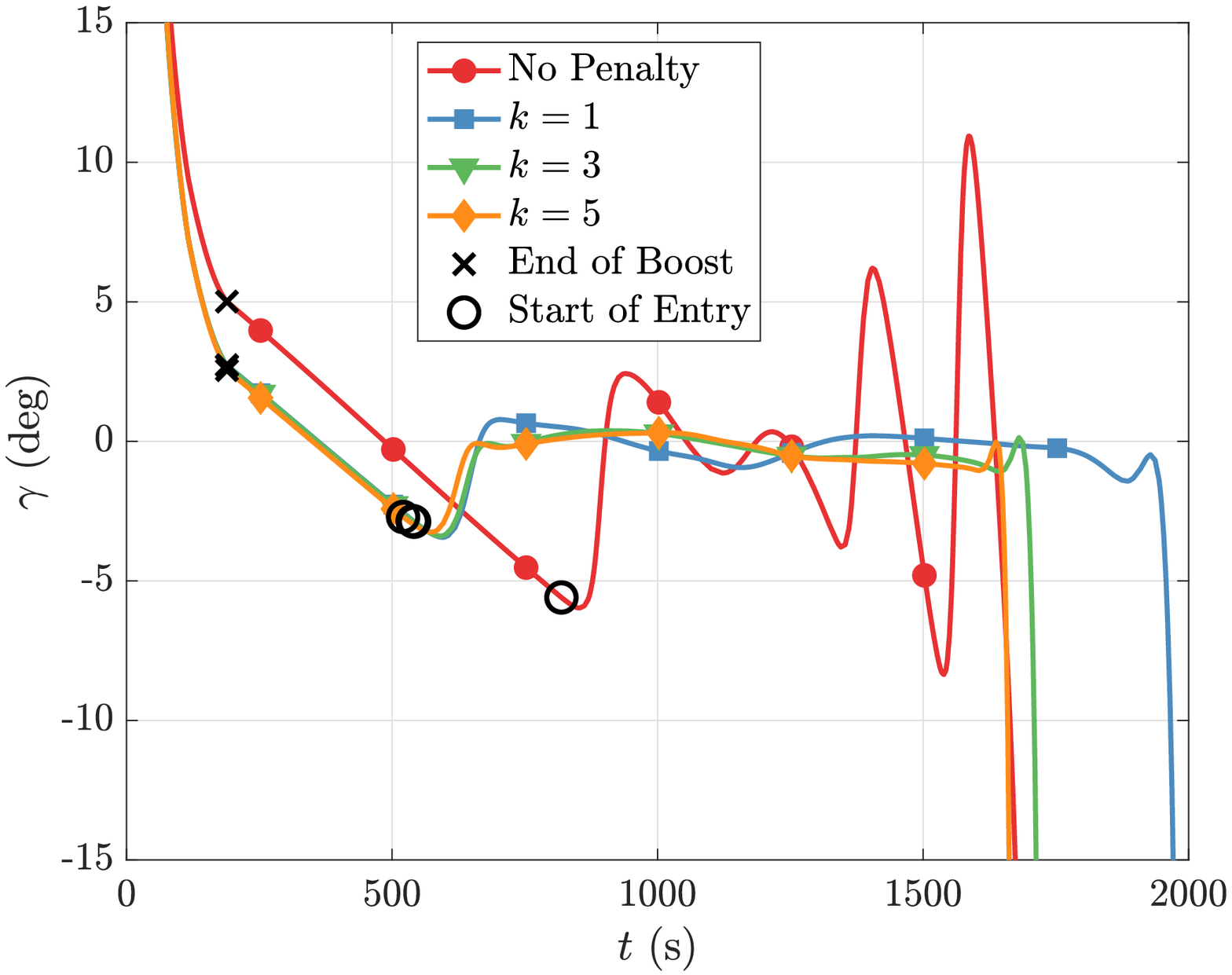}}
 \end{tabular}
  \caption{Variations in flight path angle as $k$ is varied.}\label{fig:sig-fpa}
\end{figure}

\begin{figure}[hbt!]
  \centering
  \begin{tabular}{lr}
  \subfloat[Altitude, $h(t)$ vs. time, $t$.]{\includegraphics[width=.475\textwidth]{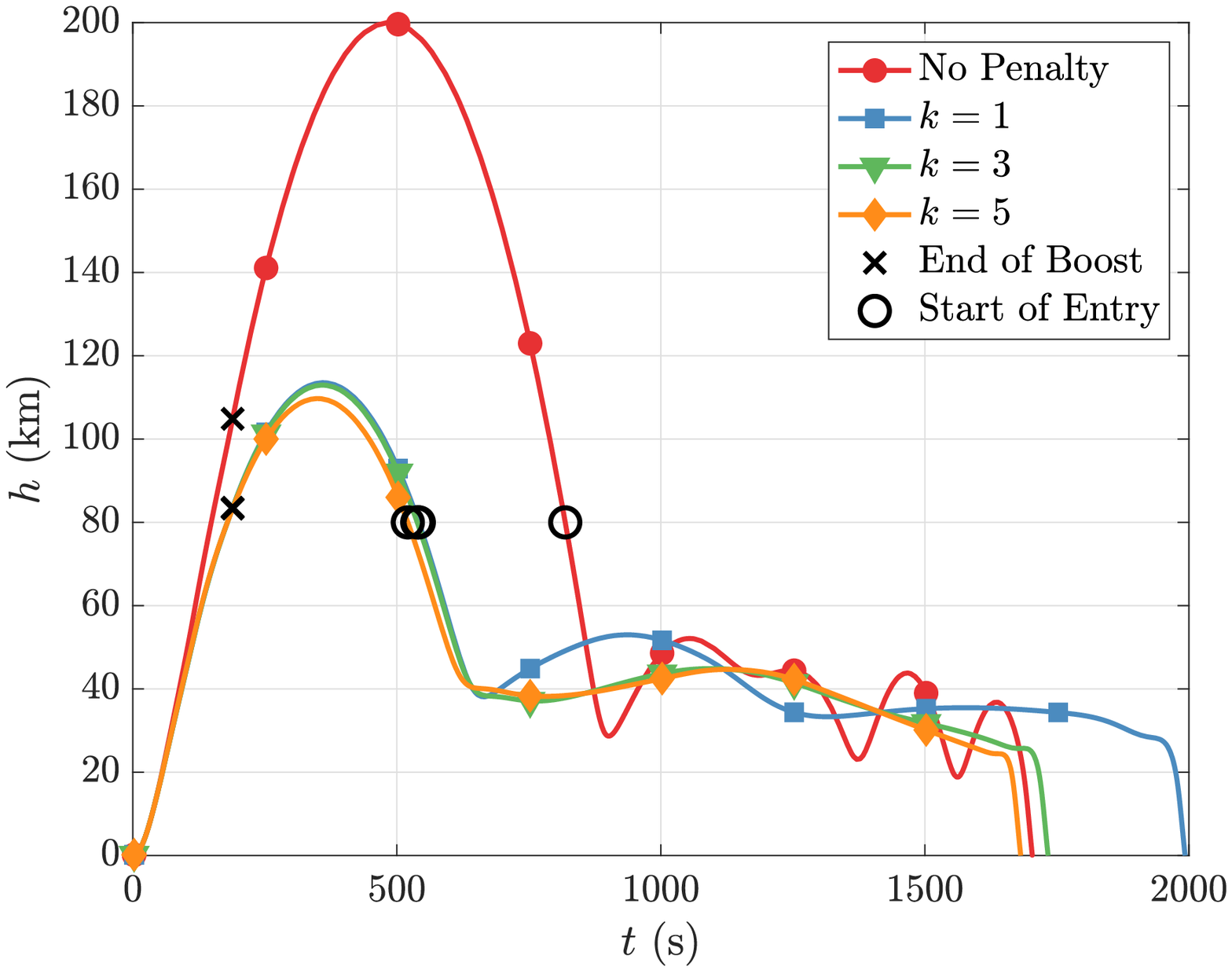}}
  &
  \subfloat[Speed, $v(t)$ vs. time, $t$.]{\includegraphics[width=.475\textwidth]{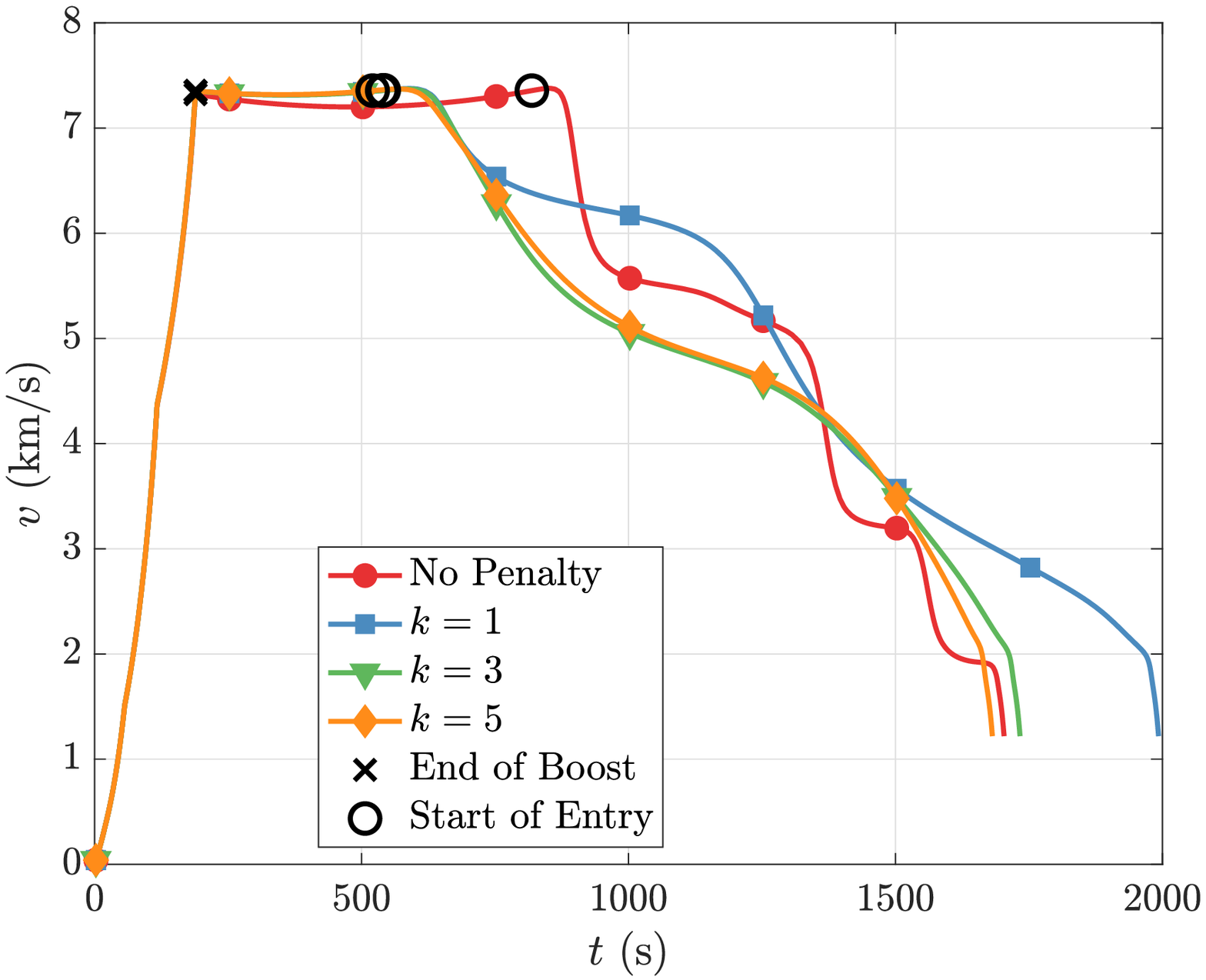}}
 \end{tabular}
  \caption{Variations in altitude and speed as $k$ is varied.}\label{fig:sig-hv}
\end{figure}

Now consider the trajectories obtained with and without the penalty term added to the integrand of Eq.~\eqref{eq:Cost-phase} in phase 8.  The flight path angle profiles of each nominal trajectory are shown in Fig.~\ref{fig:sig-fpa}.  Observing Fig.~\ref{fig:sig-fpa} it is seen that the addition of the penalty term produces trajectories with flight path angles much closer to zero during entry when compared to the large oscillations in the flight path angle observed for the solution obtained without the penalty term.  The increase in glide-like behavior with the addition of the penalty term is further evidenced in Fig.~\ref{fig:sig-hv} where it is seen that the phugoidal oscillations in altitude observed for the no penalty solution are largely gone in the $k = \{1, 3, 5\}$ solutions.  Figure~\ref{fig:sig-hv} also shows that speed is more smoothly depleted during entry for the $k = \{1, 3, 5\}$ solutions as opposed to the step-like reductions in speed seen in the no penalty solution.  Interestingly, in Figs.~\ref{fig:sig-fpa} and \ref{fig:sig-hv} the $k = \{1, 3, 5\}$ solutions are all qualitatively similar to one another and the $k = 3$ and $k = 5$ solutions are nearly identical.

\begin{figure}[hbt!]
  \centering
  \begin{tabular}{lr}
  \subfloat[Heating rate, $\dot{Q}(t)$ vs. time, $t$.]{\includegraphics[width=.475\textwidth]{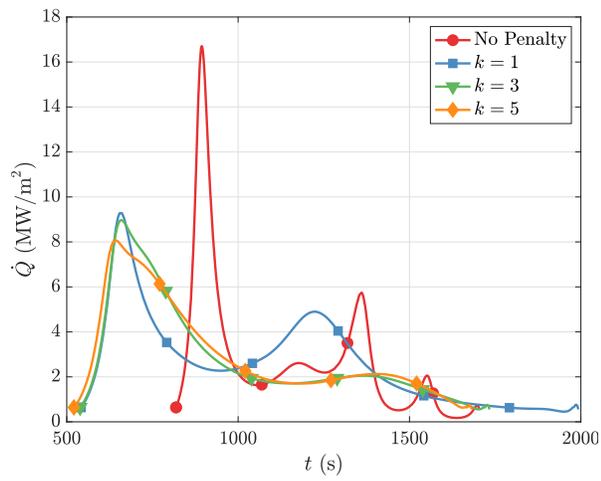}}
  &
  \subfloat[Sensed acceleration, $n(t)$ vs. time, $t$.]{\includegraphics[width=.475\textwidth]{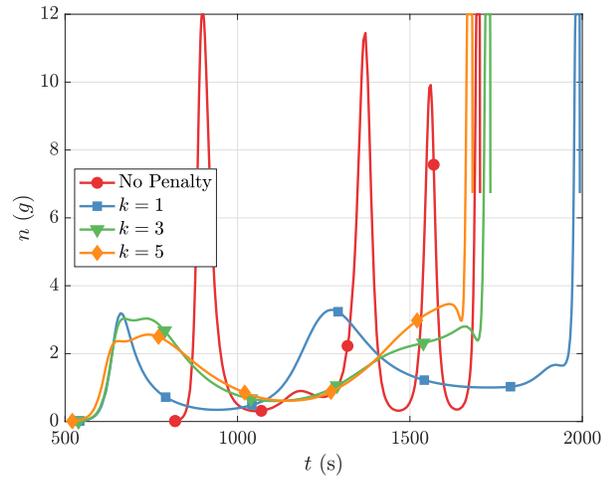}}
 \end{tabular}
  \caption{Variations in the heating rate and sensed acceleration profiles during entry as $k$ is varied.}\label{fig:sig-Qdotn}
\end{figure}

The phugoid oscillations present in the solution (or lack thereof) directly impact the sensed acceleration, dynamic pressure, and heating rate experienced by the entry vehicle.  Figure~\ref{fig:sig-Qdotn} compares the changes in constraint behavior for the heating rate and sensed acceleration profiles, noting that the dynamic pressure and sensed acceleration profiles observed during entry are qualitatively similar to one another.  It is apparent in Fig.~\ref{fig:sig-Qdotn} that the no penalty solution sees large spikes in the heating rate and sensed acceleration relative to the $k = \{1, 3, 5\}$ solutions.  Comparing the no penalty solution in Figs.~\ref{fig:sig-hv} and \ref{fig:sig-Qdotn}, it is noticed that the spikes in heating rate and sensed acceleration coincide with each descent into the Earth's atmosphere where speed is reduced rapidly.  In contrast, the more gradual changes in altitude and speed observed in Fig.~\ref{fig:sig-hv} for the $k = \{1, 3, 5\}$ solutions coincides with the more evenly distributed heating rate and sensed acceleration profiles observed in Fig.~\ref{fig:sig-Qdotn}.  Thus, the phugoid penalty term accomplishes the goal of smoothing out the entry trajectory so that the heating rate, sensed acceleration, and dynamic pressure experienced by the vehicle during entry is more evenly distributed in time.

\begin{figure}[hbt!]
  \centering
  \begin{tabular}{lr}
  \subfloat[Angle of attack, $\alpha(t)$ vs. time, $t$.]{\includegraphics[width=.475\textwidth]{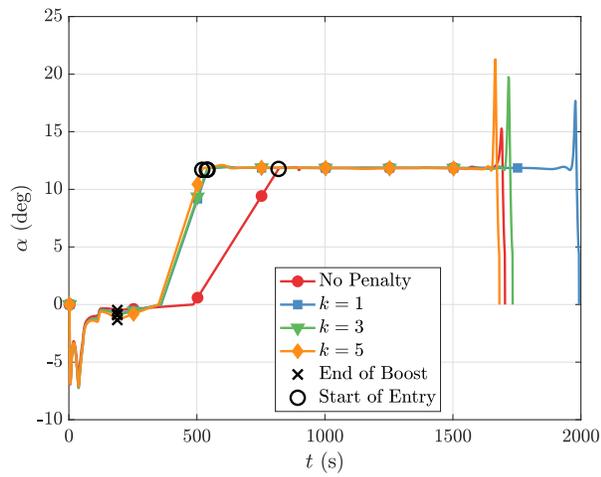}}
  &
  \subfloat[Bank angle, $\sigma(t)$ vs. time, $t$.]{\includegraphics[width=.475\textwidth]{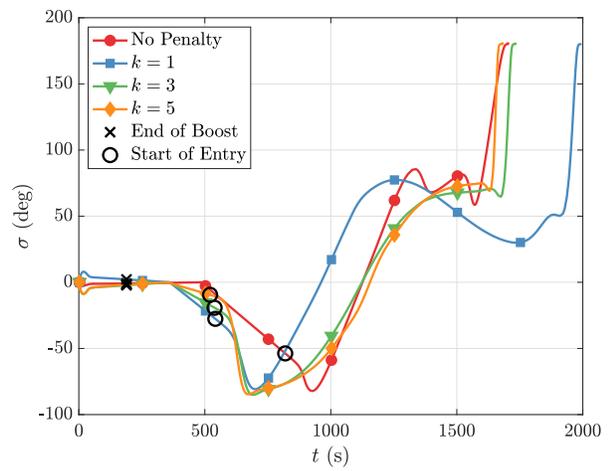}}
 \end{tabular}
  \caption{Optimal angle of attack and bank angle profiles obtained as $k$ is varied.}\label{fig:sig-aoabank}
\end{figure}

While the previous results motivate the use of the phugoid penalty term defined in Eq.~\eqref{eq:sig}, it is important to determine whether the addition of the penalty term adversely affects the goal of generating trajectories with wide control margins.  Figure~\ref{fig:sig-aoabank} shows the angle of attack and bank angle profiles obtained for each of the solutions, noting that Eq.~\eqref{eq:convert-state} has been used to convert the Euler parameters to the bank angle in phases 1 and 8.  In Fig.~\ref{fig:sig-aoabank} it is seen that the angle of attack and bank angle profiles are all qualitatively similar to one another, with the most rapid changes occurring at the end of entry as the vehicle turns over and dives onto the target.  While the peak angle of attack required for the turnover maneuver does increase with the addition of the penalty term (and as $k$ increases), the overall performance remains similar to that of the no penalty solution.

%
%
%
%
%

\clearpage	
\subsection{Constrained Heating Rate and Heating Load Studies}\label{sect:conStudy}
The nominal solution behavior of the ascent-entry trajectory optimization problem stated in Section~\ref{sect:ProblemFormulation} is studied in Section~\ref{sect:sigStudy}.  In particular, Section~\ref{sect:sigStudy} explores the effects of adding a phugoid oscillation penalty term to the cost while assuming unconstrained heating rate and heating load limits.  This section builds upon Section~\ref{sect:sigStudy} with three separate studies, each employing the phugoid penalty term with $k = 3$, and each characterizing changes in nominal performance under varying maximum heating rate and heating load requirements.  First, Section~\ref{sect:performance-HR} explores variations in performance with $Q_{\max} = \infty$ and as $\dot{Q}_{\max}$ is lowered by $1~\textrm{MW/m}^2$ increments from the nominal solution until failure (the trajectory optimization procedure fails to converge to a solution).  Next, Section~\ref{sect:performance-Q} studies the performance with $\dot{Q}_{\max} = \infty$ and as $Q_{\max}$ is lowered by $100~\textrm{MJ/m}^2$ increments from the nominal solution until failure.  Finally, Section~\ref{sect:performance-combined} compares the behavior of three combined constraint ($\dot{Q}_{\max}$ and $Q_{\max}$ both constrained) missions of varying difficulty.

\subsubsection{Heating Load Unconstrained Solutions}\label{sect:performance-HR}
\begin{figure}[hbt!]
  \centering
  \includegraphics[width=.475\textwidth]{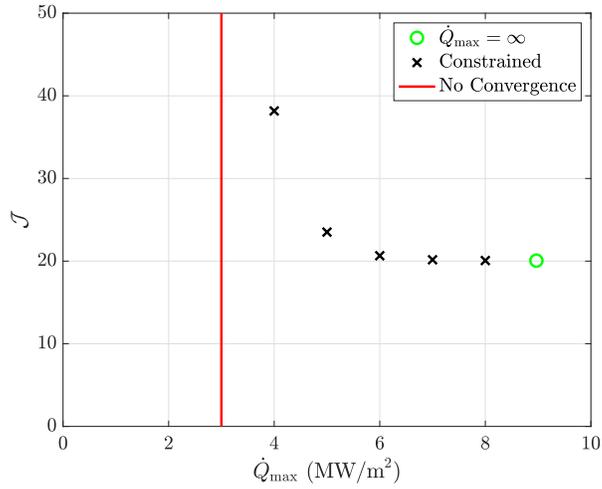}
  \caption{Performance as $\dot{Q}_{\max}$ is lowered.}\label{fig:Performance-HRcost}
\end{figure}

\begin{figure}[hbt!]
  \centering
  \includegraphics[width=.475\textwidth]{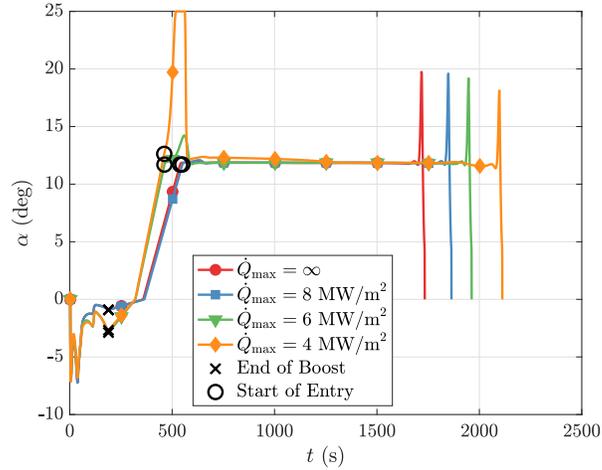}
  \caption{Optimal angle of attack profiles as $\dot{Q}_{\max}$ is lowered.}\label{fig:Performance-HRaoa}
\end{figure}

Figure~\ref{fig:Performance-HRcost} shows the change in value of the performance index as a function of the maximum allowable stagnation point heating rate.  The nominal solution attains a maximum heating rate of $9.0~\textrm{MW}/\textrm{m}^2$ and failure occurs at $3~\textrm{MW}/\textrm{m}^2$.  It is also observed that solutions above $6~\textrm{MW}/\textrm{m}^2$ attain similar values of the performance index as the nominal solution.  This similarity in performance is further evidenced by Fig.~\ref{fig:Performance-HRaoa} where it is seen that the angle of attack profiles for $\dot{Q}_{\max} = \{6, 8, \infty\}~\textrm{MW}/\textrm{m}^2$ are quite similar to each other in both the boost and glide phases.  In contrast, the solution obtained at $\dot{Q}_{\max} = 4~\textrm{MW}/\textrm{m}^2$ requires the angle of attack to be at its peak near the start of entry.

\begin{figure}[hbt!]
  \centering
  \begin{tabular}{lr}
  \subfloat[Altitude, $h(t)$ vs. time, $t$.]{\includegraphics[width=.475\textwidth]{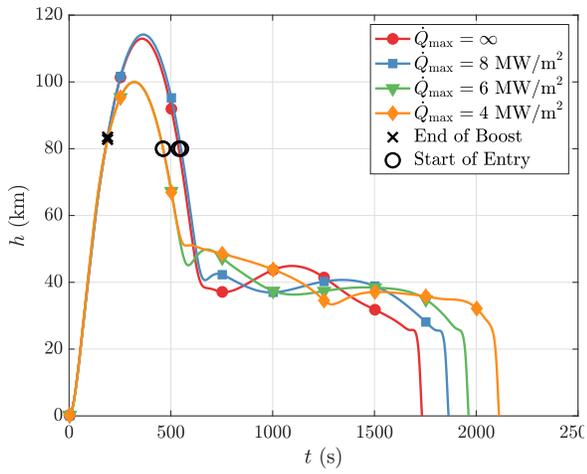}}
  &
  \subfloat[Speed, $v(t)$ vs. time, $t$.]{\includegraphics[width=.475\textwidth]{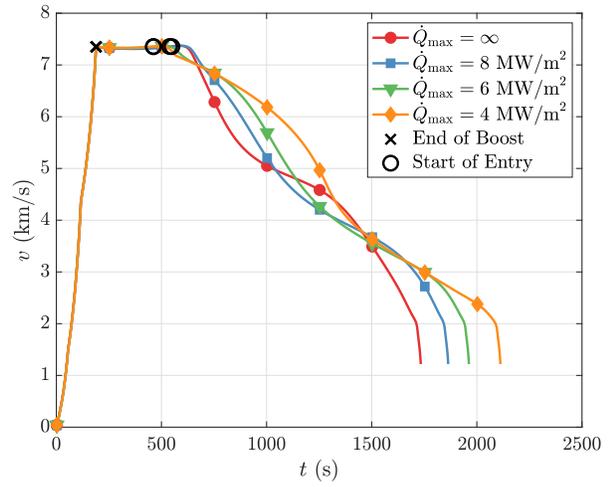}}
 \end{tabular}
  \caption{Variations in altitude and speed as $\dot{Q}_{\max}$ is lowered.}\label{fig:Performance-HRhv}
\end{figure}

It is interesting to note, however, that while the value of the performance index is similar above $6~\textrm{MW}/\textrm{m}^2$, the corresponding trajectories do exhibit different behavior from one another.  Figure~\ref{fig:Performance-HRhv} shows the optimal altitude and speed profiles attained for $\dot{Q}_{\max} = \{4, 6, 8, \infty\}~\textrm{MW}/\textrm{m}^2$.  It is observed that the peak altitude, $h\left(t_0^{(6)}\right)$, decreases to its lower bound ($100~\textrm{km}$) as $\dot{Q}_{\max}$ is lowered.  The lower peak altitude allows for a more shallow entry flight path angle as evidenced by the data in Table~\ref{tab:Performance-HR}.  The smaller flight path angle, in turn, allows for the entry vehicle to maintain a higher altitude upon initial descent into the atmosphere.  Maintaining sufficiently high altitude on initial descent is critical because $\dot{Q}$ is a function of both altitude and speed, and speed (roughly $7.35~\textrm{km}/\textrm{s}$ at the pierce point) cannot be rapidly reduced until sufficiently low altitudes are reached.

\begin{figure}[hbt!]
  \centering
  \includegraphics[width=.475\textwidth]{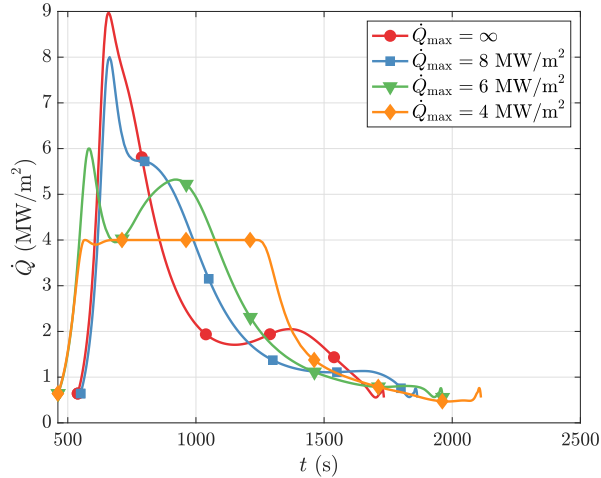}
  \caption{Stagnation point heating rate profiles during entry as $\dot{Q}_{\max}$ is lowered.}\label{fig:Performance-HRQdot}
\end{figure}

Next, Fig.~\ref{fig:Performance-HRQdot} shows the structure of the stagnation point heating rate profiles attained for $\dot{Q}_{\max} = \{4, 6, 8, \infty\}~\textrm{MW}/\textrm{m}^2$.  It is observed that the $\dot{Q}_{\max}$ limit is reached at a single touch-and-go point for $\dot{Q}_{\max} = \{6, 8, \infty\}~\textrm{MW}/\textrm{m}^2$, whereas the solution obtained for $\dot{Q}_{\max} = 4~\textrm{MW}/\textrm{m}^2$ exhibits a touch-and-go point followed by an arc of the solution along $\dot{Q}_{\max}$.  Observing Fig.~\ref{fig:Performance-HRhv} once more, it is seen that each touch-and-go point occurs at or near the lowest altitude achieved during initial descent, as the entry vehicle transitions from ballistic descent to atmospheric flight.  It is also observed that the $\dot{Q}_{\max} = 4~\textrm{MW}/\textrm{m}^2$ solution exhibits monotonically decreasing altitude during the $\dot{Q} = \dot{Q}_{\max}$ arc.

Finally, it is noted that while lower values of $\dot{Q}_{\max}$ might be thought to decrease the total heating load $Q$ as well, the opposite trend is observed in Table~\ref{tab:Performance-HR}.  Upon further examination of Table~\ref{tab:Performance-HR} it is seen that the total time elapsed during entry grows with decreasing $\dot{Q}_{\max}$ as well.  Thus, one contributor to the increasing heat load is simply the increased time during entry over which the heat load accumulates.  It is also noted that while the peak value of $\dot{Q}$ decreases with lower $\dot{Q}_{\max}$ values, it is not necessarily the case that the average value of $\dot{Q}$ decreases as well.  However, the trajectories obtained at $\dot{Q}_{\max} = \{\infty, 8, 7, 6, 5, 4\}~\textrm{MW}/\textrm{m}^2$ did exhibit simultaneously decreasing $\dot{Q}_{\max}$ and average $\dot{Q}$ values during entry, noting that the respective averages are $\{ 2.85, 2.73, 2.72, 2.62, 2.49, 2.40 \}~\textrm{MW}/\textrm{m}^2$.

\begin{table}[h]
  \centering
  \caption{Numerical results for the unconstrained heating load solutions\label{tab:Performance-HR}}
  \renewcommand{\baselinestretch}{1}\normalsize\normalfont
 \begin{tabular}{ccccccc}\hline
     $\dot{Q}_{\max}$ & 
     $Q$ & 
     $\C{J}$ & 
     $h\left(t_0^{(6)}\right)$ & 
     $v\left(t_0^{(7)}\right)$ & 
     $\gamma\left(t_0^{(7)}\right)$ & 
     $t_f^{(8)} - t_0^{(7)}$ \\
     
     $\left(\textrm{MW/m}^2\right)$ & 
     $\left(\textrm{MJ/m}^2\right)$ & 
      & 
      $\left(\textrm{km}\right)$ & 
      $\left(\textrm{km/s}\right)$ & 
      $\left(\textrm{deg}\right)$ & 
      $\left(\textrm{s}\right)$ \\\hline
     
     $\infty$ &
      3400 & 20.06 & 113.0 & 7.36 & -2.85 & 1195 \\
     
     $8$ & 
     3581 & 20.07 & 114.2 & 7.35 & -2.87 & 1314 \\
     
      $7$ & 
     3809 & 20.16 & 107.8 & 7.35 & -2.59 & 1402 \\
     
      $6$ & 
     3930 & 20.65 & 100.0 & 7.35 & -2.20 & 1500 \\
     
      $5$ & 
     3896 & 23.51 & 100.0 & 7.35 & -2.20 & 1566 \\
     
     $4$ & 
     3977 & 38.17 & 100.0 & 7.36 & -2.20 & 1651 \\\hline
  \end{tabular}
\end{table}

\subsubsection{Heating Rate Unconstrained Solutions}\label{sect:performance-Q}
\begin{figure}[hbt!]
  \centering
  \includegraphics[width=.475\textwidth]{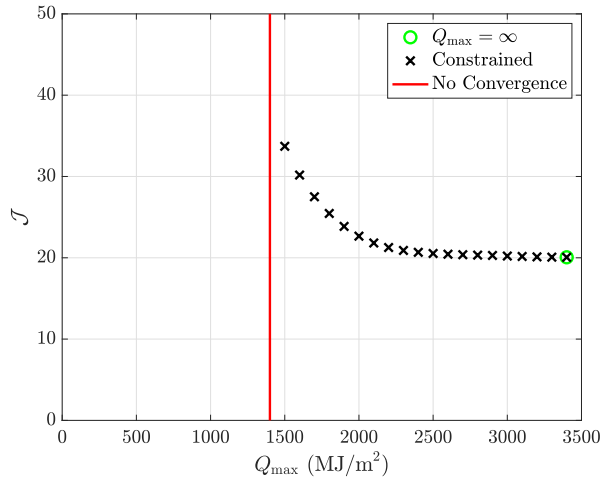}
  \caption{Performance as $Q_{\max}$ is lowered.}\label{fig:Performance-Qcost}
\end{figure}

\begin{figure}[hbt!]
  \centering
  \includegraphics[width=.475\textwidth]{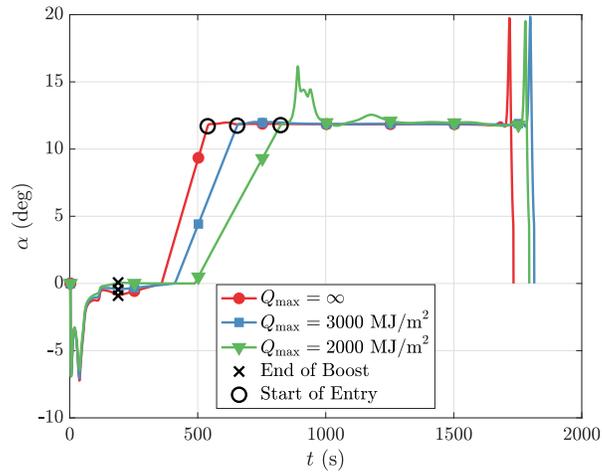}
  \caption{Optimal angle of attack profiles as $Q_{\max}$ is lowered.}\label{fig:Performance-Qaoa}
\end{figure}

Figure~\ref{fig:Performance-Qcost} shows the change in value of the performance index as a function of the maximum allowable stagnation point heating load.  The nominal solution attains a maximum heating load of $3400~\textrm{MJ}/\textrm{m}^2$ and failure occurs at $1400~\textrm{MJ}/\textrm{m}^2$.  Upon further inspection of Fig.~\ref{fig:Performance-Qcost} it is noticed that solutions above $2500~\textrm{MJ}/\textrm{m}^2$ attain similar values of the performance index as the nominal solution.  Figure~\ref{fig:Performance-Qaoa} further attests to the performance similarities, noting that the angle of attack profiles obtained for ${Q}_{\max} = \{3000, \infty\}~\textrm{MJ}/\textrm{m}^2$ are quite similar to one another, whereas the solution obtained at ${Q}_{\max} = 2000~\textrm{MJ}/\textrm{m}^2$ sees much larger fluctuations in the angle of attack at the start of entry.

\begin{figure}[hbt!]
  \centering
  \begin{tabular}{lr}
  \subfloat[Altitude, $h$ vs. time, $t$.]{\includegraphics[width=.475\textwidth]{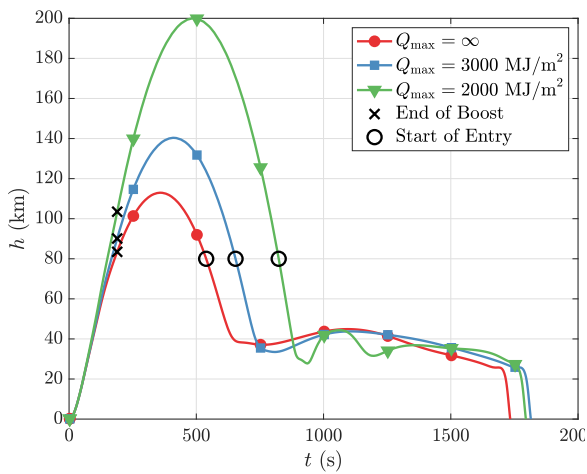}}
  &
  \subfloat[Speed, $v$ vs. time, $t$.]{\includegraphics[width=.475\textwidth]{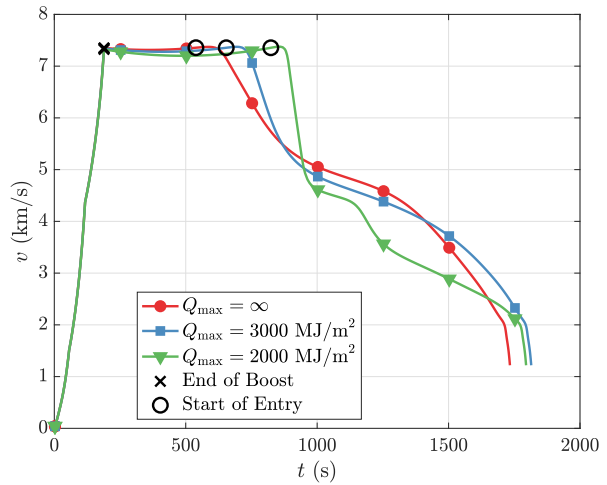}}
 \end{tabular}
  \caption{Variations in altitude and speed as $Q_{\max}$ is lowered.}\label{fig:Performance-Qhv}
\end{figure}

Next, Fig.~\ref{fig:Performance-Qhv} shows the variations in the optimal altitude and speed profiles as $Q_{\max}$ is lowered.  Inspecting Fig.~\ref{fig:Performance-Qhv} it is seen that the peak altitude $h\left(t_0^{(6)}\right)$ increases with decreasing heating load limits.  In fact, Table~\ref{tab:Performance-Q} indicates that the peak altitude is constrained to the upper limit of $200~\textrm{km}$ once $Q_{\max}$ is lowered at or below $2500~\textrm{MJ}/\textrm{m}^2$.  Another trend observed in Fig.~\ref{fig:Performance-Qhv} is that decreasing $Q_{\max}$ leads to more rapid speed depletion for the entry vehicle upon initial descent into the Earth's atmosphere.  The large reduction in speed near the beginning of entry takes advantage of the fact that $Q$ is the integral of $\dot{Q}$ which is, in turn, a function of $v^{3.15}$.  Thus, rapid speed depletion at the start of entry can enable lower average values of $\dot{Q}$ achieved during the entirety of entry even if, as Table~\ref{tab:Performance-Q} indicates, the maximum heating rate experienced during entry increases.

\begin{table}[h]
  \centering
  \caption{Numerical results for the unconstrained heating rate solutions\label{tab:Performance-Q}}
  \renewcommand{\baselinestretch}{1}\normalsize\normalfont
 \begin{tabular}{ccccccc}\hline
     $Q_{\max}$ & 
     $\max\left(\dot{Q}\right)$ & 
     $\C{J}$ & 
     $h\left(t_0^{(6)}\right)$ & 
     $v\left(t_0^{(7)}\right)$ & 
     $\gamma\left(t_0^{(7)}\right)$ & 
     $t_f^{(8)} - t_0^{(7)}$ \\
     
     $\left(\textrm{MJ/m}^2\right)$ & 
     $\left(\textrm{MW/m}^2\right)$ & 
      & 
      $\left(\textrm{km}\right)$ & 
      $\left(\textrm{km/s}\right)$ & 
      $\left(\textrm{deg}\right)$ & 
      $\left(\textrm{s}\right)$ \\\hline
     
     $\infty$ &
      9.0 & 20.06 & 113.0 & 7.36 & -2.85 & 1195 \\
     
     $3000$ & 
     11.3 & 20.21 & 140.4 & 7.36 & -3.86 & 1160 \\
     
     $2500$ & 
     14.8 & 20.53 & 200.0 & 7.35 & -5.52 & 1039 \\
     
     $2000$ & 
     14.6 & 22.66 & 200.0 & 7.35 & -5.57 & 971 \\
     
     $1500$ & 
     11.8 & 33.7 & 200.0 & 7.36 & -5.59 & 887 \\\hline
  \end{tabular}
\end{table}

\newpage 
\subsubsection{Combined-Constraint Missions}\label{sect:performance-combined}
Sections~\ref{sect:performance-HR} and \ref{sect:performance-Q} analyze changes in the combined ascent-entry solution behavior as $\dot{Q}_{\max}$ is varied individually while $Q_{\max} = \infty$ or vice-versa.  Some of the trends observed in Sections~\ref{sect:performance-HR} and \ref{sect:performance-Q}, such as increases in heating load as $\dot{Q}_{\max}$ is lowered, indicate that performance trade-offs exist when enforcing both the heating rate and heating load constraints simultaneously.  This section endeavors to characterize such trade-offs, as well as to analyze the solution behavior of the combined-constraint mission to see if any of the trends observed in Sections~\ref{sect:performance-HR} and \ref{sect:performance-Q} reappear or if new trends appear.

\begin{figure}[hbt!]
  \centering
  \includegraphics[width=.475\textwidth]{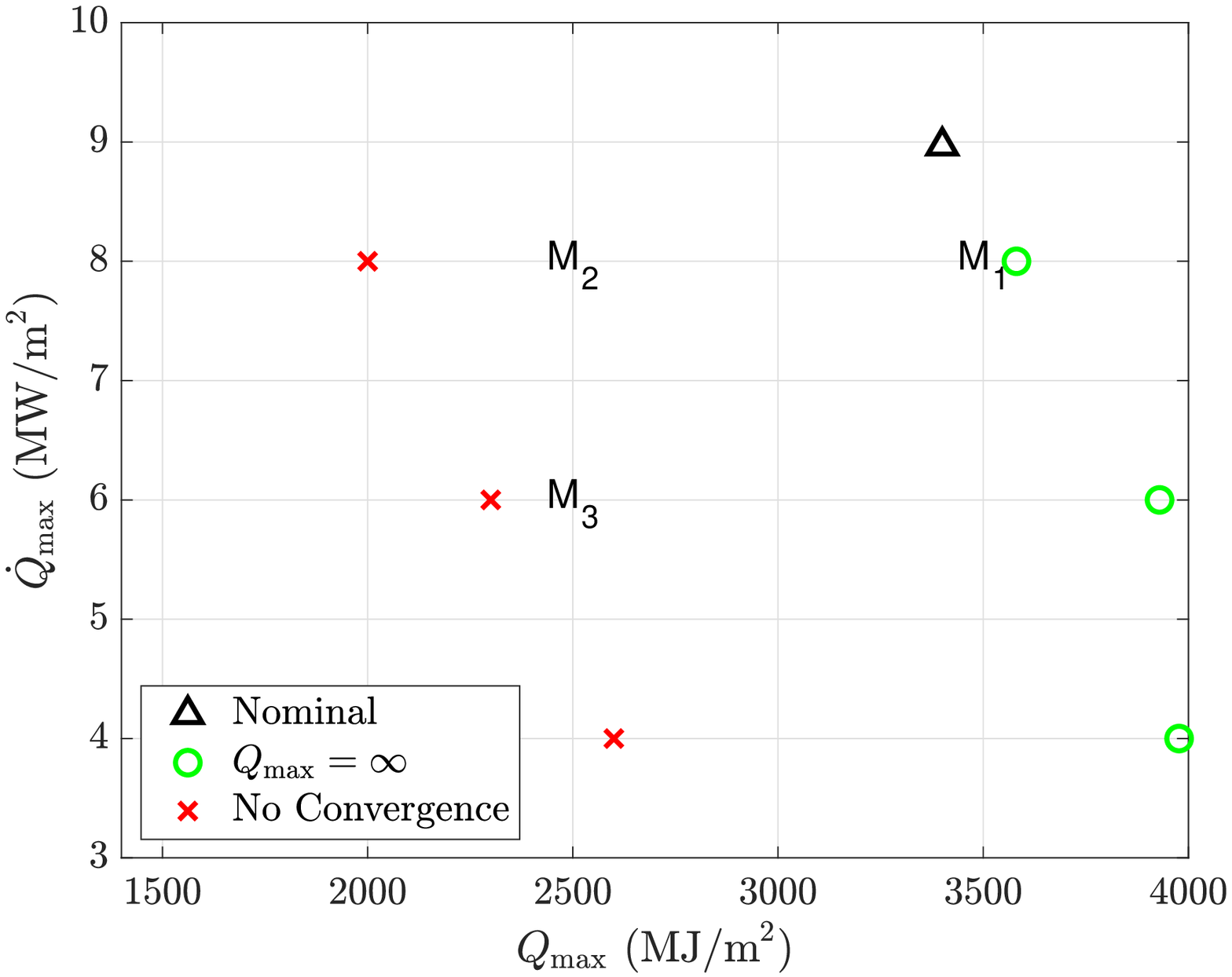}
  \caption{Failure points for combined constraint solutions relative to the nominal and $Q_{\max} = \infty$ solutions.  Note that $\{\textrm{M}_1, \textrm{M}_2, \textrm{M}_3\}$ denote Missions 1-3 respectively.}\label{fig:FailurePoints}
\end{figure}

Figure~\ref{fig:FailurePoints} provides a sense of the trade-offs in performance for combined-constraint missions.  The nominal and $Q_{\max} = \infty$ solutions act as reference points in the heating rate and heating load constraint space, and each of the failure points are obtained by lowering $Q_{\max}$ by $100~\textrm{MJ}/\textrm{m}^2$ increments (while holding $\dot{Q}_{\max}$ constant) until the trajectory optimization routine fails to converge to a solution.  Clearly, lower $\dot{Q}_{\max}$ limits come at the expense of raising the smallest feasible $Q_{\max}$ value (and vice-versa).  It is also observed that $Q_{\max}$ may be lowered significantly from the natural heating loads obtained when the heating load is unconstrained ($Q_{\max} = \infty$).  Taken together, Fig.~\ref{fig:FailurePoints} provides a useful map for designing combined-constraint mission requirements and for gauging the relative difficulty of a given mission.

\begin{figure}[hbt!]
  \centering
  \begin{tabular}{lr}
  \subfloat[Heating Rate, $\dot{Q}(t)$ vs. time, $t$.]{\includegraphics[width=.475\textwidth]{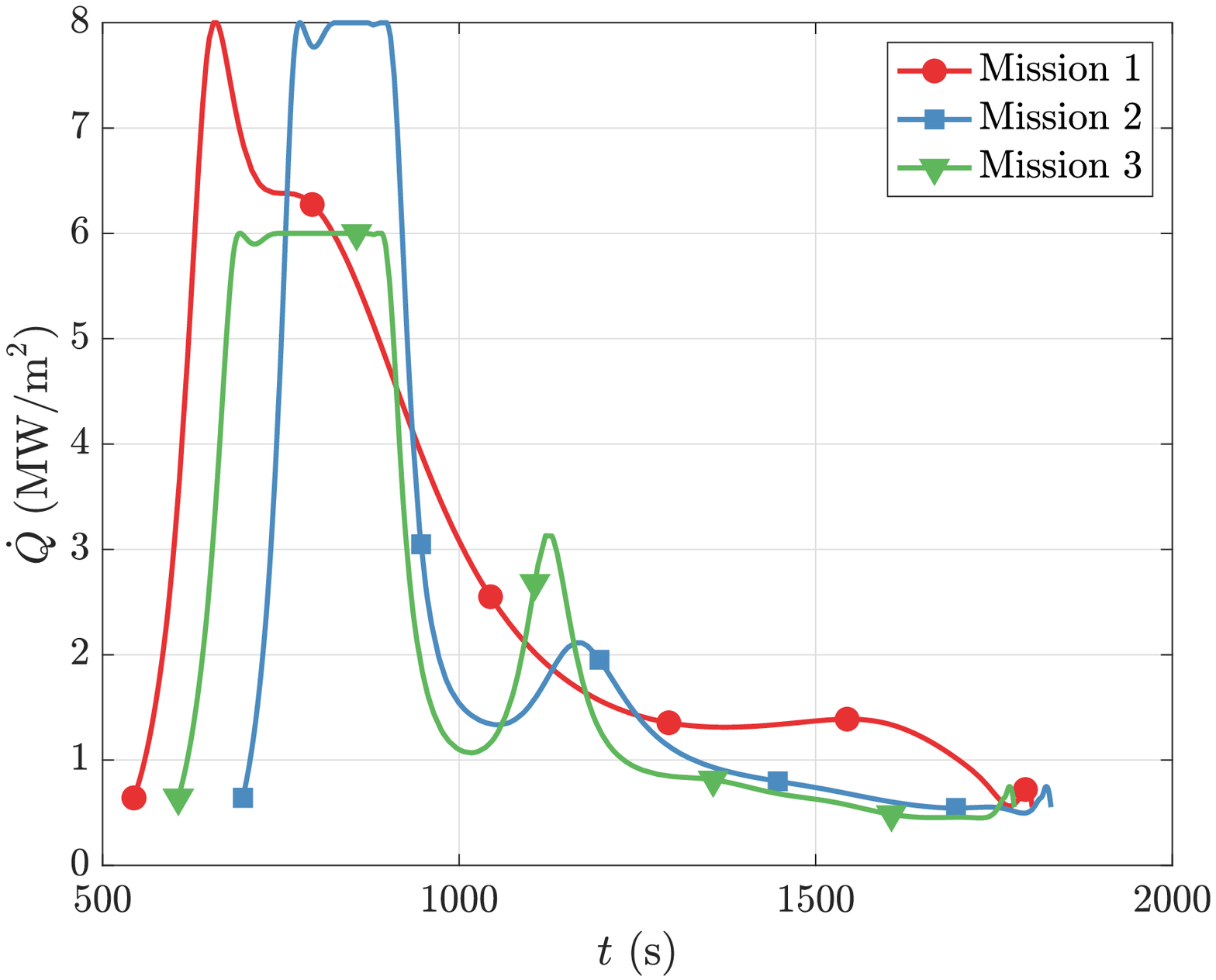}}
  &
  \subfloat[Sensed Acceleration, $n(t)$ vs. time, $t$.]{\includegraphics[width=.475\textwidth]{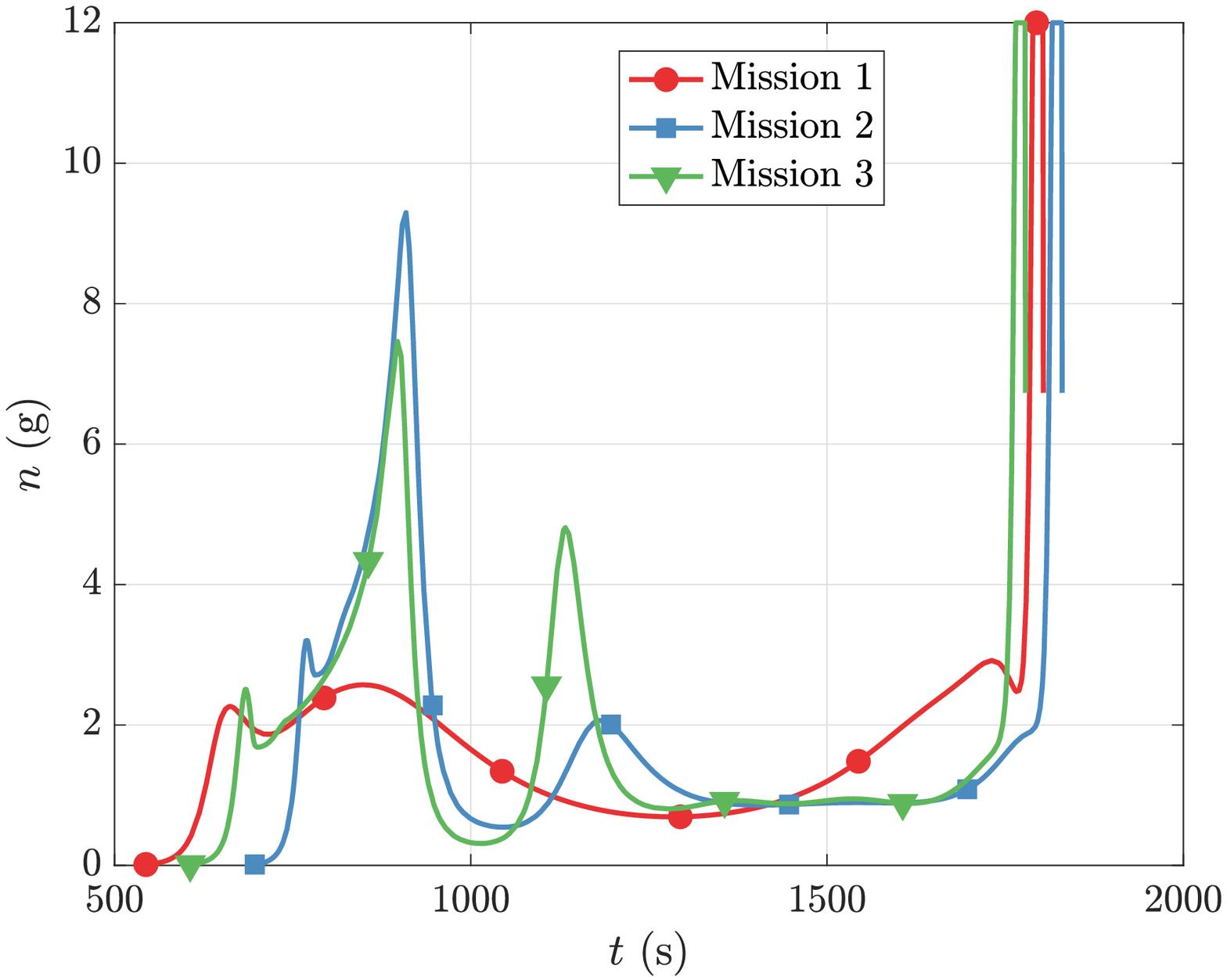}}
 \end{tabular}
  \caption{Heating rate and sensed acceleration profiles obtained during entry for Missions 1-3.}\label{fig:Combined-Qdotn}
\end{figure}

Next, consider the following three combined constraint missions of easy, medium, and hard difficulty (as measured by the performance index and their proximity to the failure points shown in Fig.~\ref{fig:FailurePoints}).  The three missions are listed in order of increasing difficulty as
\begin{itemize}
	\item \textbf{Mission 1:}~~~$\{ \dot{Q}_{\max} = 8~\textrm{MW}/\textrm{m}^2~,~{Q}_{\max} = 3500~\textrm{MJ}/\textrm{m}^2\}$,
	\item \textbf{Mission 2:}~~~$\{ \dot{Q}_{\max} = 8~\textrm{MW}/\textrm{m}^2~,~{Q}_{\max} = 2500~\textrm{MJ}/\textrm{m}^2\}$,
	\item \textbf{Mission 3:}~~~$\{ \dot{Q}_{\max} = 6~\textrm{MW}/\textrm{m}^2~,~{Q}_{\max} = 2500~\textrm{MJ}/\textrm{m}^2\}$.
\end{itemize}
First, observe in Fig.~\ref{fig:Combined-Qdotn} the points/arcs of each solution where the heating rate and sensed acceleration constraints are active.  Similar to Section~\ref{sect:performance-HR}, it is observed that the max heating rate limit is reached at a touch-and-go point for Mission 1 (higher $\dot{Q}_{\max}$ value) and the heating rate constraint is active at a touch-and-go point followed by an arc for Mission 3 (lower $\dot{Q}_{\max}$ value).  Interestingly, Mission 2 displays the same $\dot{Q}$ structure as Mission 3, indicating that the transition from touch-and-go point to a touch-and-go plus arc sequence can occur as the $Q_{\max}$ limit is lowered as well.  It is also observed in Fig.~\ref{fig:Combined-Qdotn} that Missions 2 and 3 ($Q_{\max} = 2500~\textrm{MJ}/\textrm{m}^2$) see larger acceleration loads during the first half of entry when compared to Mission 1 ($Q_{\max} = 3500~\textrm{MJ}/\textrm{m}^2$).  Finally, in Fig.~\ref{fig:Combined-Qdotn} it is seen that all three missions reach the $12~\textrm{g}$ sensed acceleration limit along a short arc as the entry vehicle carries out a dive maneuver prior to target impact.

\begin{figure}[hbt!]
  \centering
  {\includegraphics[width=.475\textwidth]{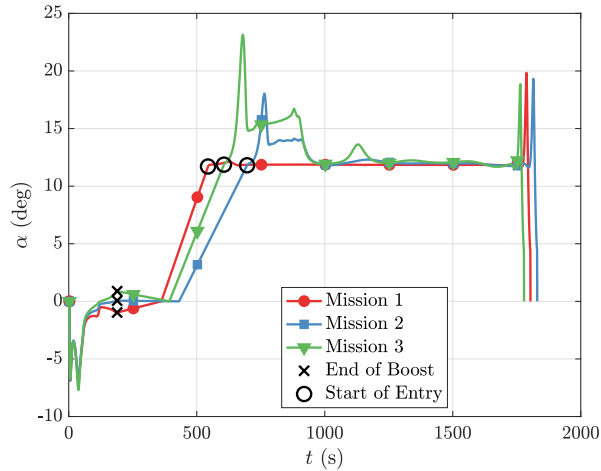}}
  \caption{Optimal angle of attack profiles obtained for Missions 1-3.}\label{fig:Combined-aoa}
\end{figure}

Next, consider the angle of attack profiles obtained for each mission and shown in Fig.~\ref{fig:Combined-aoa}.  Clearly, as the heating rate and heating load constraints tighten, the angle of attack maneuvers become more pronounced during the entry phases.  In particular, notice that the initial rise in angle of attack at the start of entry increases dramatically from Mission 1 to Mission 3.  It is also interesting to note that the angle of attack during the boost phases remains similar across all three missions, noting that the largest differences occur near the end of the boost phases.

\begin{figure}[hbt!]
  \centering
  \begin{tabular}{lr}
  \subfloat[Altitude, $h$ vs. time, $t$.]{\includegraphics[width=.475\textwidth]{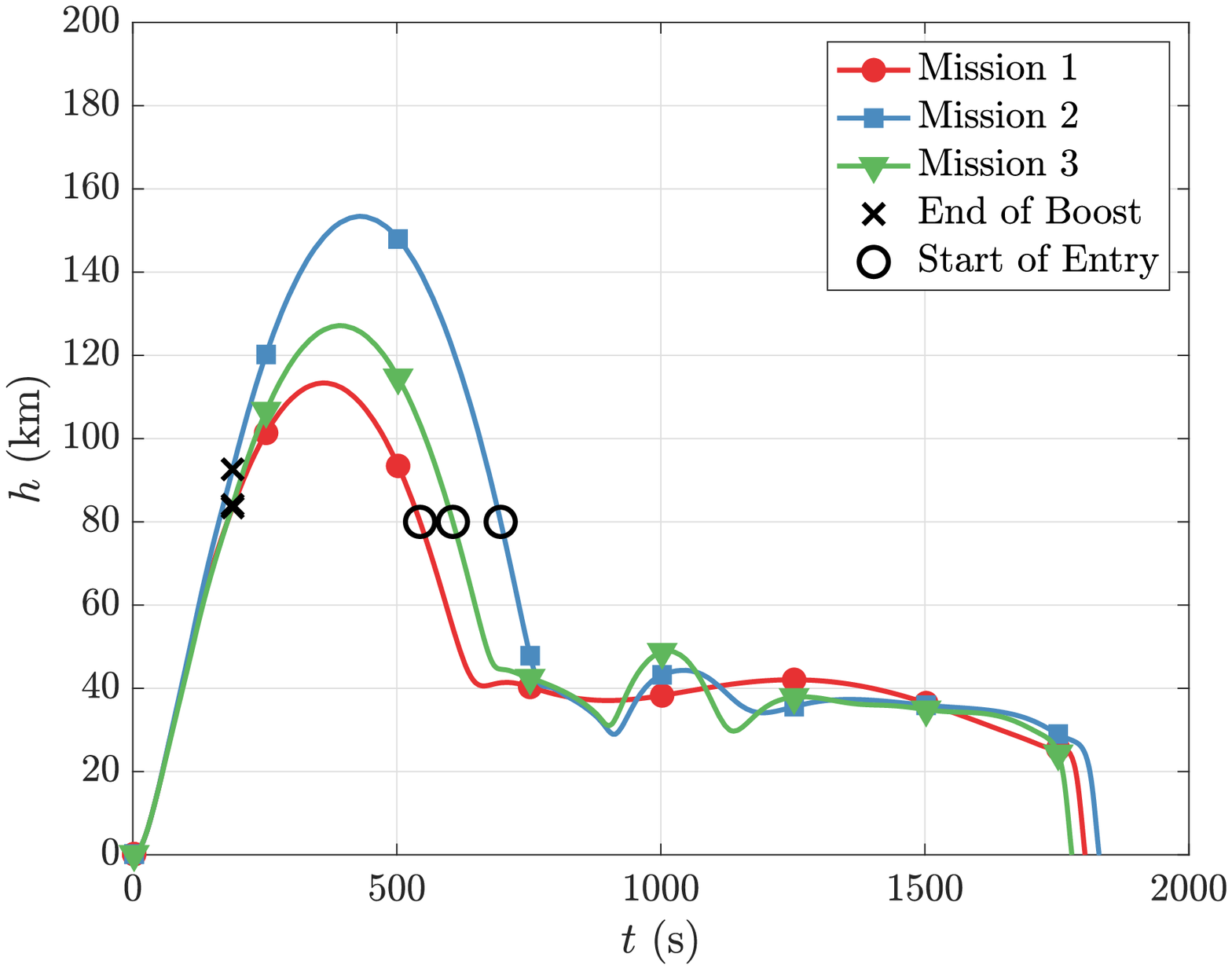}}
  &
  \subfloat[Speed, $v$ vs. time, $t$.]{\includegraphics[width=.475\textwidth]{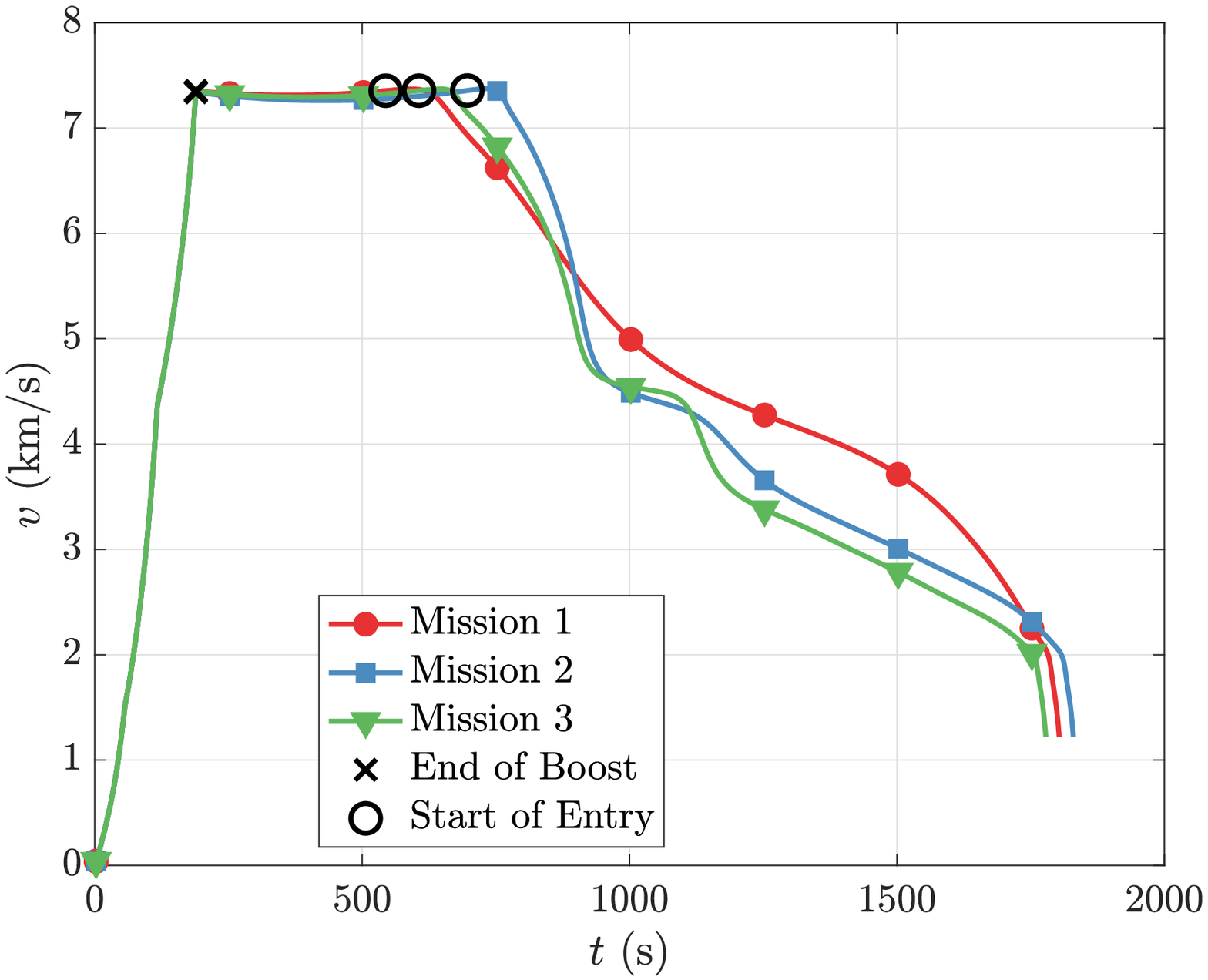}}
 \end{tabular}
  \caption{Altitude and speed profiles obtained for Missions 1-3.}\label{fig:Combined-hv}
\end{figure}

Now consider the differences in the altitude and speed profiles obtained for each mission and shown in Fig.~\ref{fig:Combined-hv}.  Similar to Section~\ref{sect:performance-Q}, the peak altitude increases from Mission 1 to Mission 2 as $Q_{\max}$ is lowered and $\dot{Q}_{\max}$ is held constant.  Likewise, the trend of decreasing peak altitude observed in Section~\ref{sect:performance-HR} is seen from Mission 2 to Mission 3 as $\dot{Q}_{\max}$ is lowered and $Q_{\max}$ is held constant.  Figure~\ref{fig:Combined-hv} also shows that Missions 2 and 3 ($Q_{\max} = 2500~\textrm{MJ}/\textrm{m}^2$) dive deeper into the Earth's atmosphere and deplete speed more rapidly during the initial descent maneuver when compared with Mission 1 ($Q_{\max} = 3500~\textrm{MJ}/\textrm{m}^2$).  However, it is important to note that the initial descent dive is blunted in Missions 2 and 3 right at the point where the heating rate constraint becomes active.

Finally, Table~\ref{tab:Combined} summarizes the key numerical results obtained for Missions 1-3.  Investigating Table~\ref{tab:Combined} reveals that several trends observed in Sections~\ref{sect:performance-HR} and \ref{sect:performance-Q} reappear as $\dot{Q}_{\max}$ and $Q_{\max}$ are independently varied across Missions 1-3.  For instance, the pierce point flight path angle, $\gamma\left(t_0^{(7)}\right)$, and the total entry time, $t_f^{(8)} - t_0^{(7)}$, both decrease as $Q_{\max}$ is lowered and $\dot{Q}_{\max}$ is held constant from Mission 1 to Mission 2.  Likewise, the opposite change occurs from Mission 2 to Mission 3 as $\dot{Q}_{\max}$ is lowered and $Q_{\max}$ is held constant.

\begin{table}[h]
  \centering
  \caption{Numerical results for combined constraint missions\label{tab:Combined}}
  \renewcommand{\baselinestretch}{1}\normalsize\normalfont
 \begin{tabular}{lccccccc}\hline
 	Description & 
     $Q_{\max}$ & 
     $\dot{Q}_{\max}$ & 
     $\C{J}$ & 
     $h\left(t_0^{(6)}\right)$ & 
     $v\left(t_0^{(7)}\right)$ & 
     $\gamma\left(t_0^{(7)}\right)$ & 
     $t_f^{(8)} - t_0^{(7)}$ \\
     
      & 
     $\left(\textrm{MJ/m}^2\right)$ & 
     $\left(\textrm{MW/m}^2\right)$ & 
      & 
      $\left(\textrm{km}\right)$ & 
      $\left(\textrm{km/s}\right)$ & 
      $\left(\textrm{deg}\right)$ & 
      $\left(\textrm{s}\right)$ \\\hline
     
     Nominal & 
     $\infty$ & $\infty$ & 
      20.06 & 113.0 & 7.36 & -2.85 & 1195 \\
     
     Mission 1 & 
     $3500$ & $8.0$ & 
     20.07 & 113.4 & 7.35 & -2.85 & 1260 \\
     
     Mission 2 & 
     $2500$ & $8.0$ & 
     23.55 & 153.4 & 7.36 & -4.29 & 1133 \\
     
     Mission 3 & 
     $2500$ & $6.0$ & 
     32.06 & 127.2 & 7.36 & -3.44 & 1172 \\\hline
  \end{tabular}
\end{table}

\section{Discussion\label{sect:discussion}}
The results of Section~\ref{sect:results} highlight several key aspects affecting the trajectory design for the ascent-entry problem stated in Section~\ref{sect:ProblemFormulation}.  First, in Section~\ref{sect:sigStudy} it was found that the nominal trajectory ($\dot{Q}_{\max} = Q_{\max} = \infty$) obtained using the performance index of Section~\ref{sect:Cost} includes phugoid oscillations during entry.  However, it was also found that the addition of the phugoid penalty term of Eq.~\eqref{eq:sig} to the Lagrange cost in phase 8 resulted in nominal trajectories being generated with more glide-like behavior during entry.  Thus, large oscillations in the altitude and flight path angle disappeared, speed depletion occurred more steadily, and large spikes in the heating rate, sensed acceleration, and dynamic pressure were reduced (particularly near the start of entry) for the trajectories obtained with the penalty term.  Finally, while the phugoid penalty design parameter $k$ did offer some control over the shape of the trajectory, the values of $k$ tested yielded qualitatively similar results.

Next, Section~\ref{sect:conStudy} explores changes in the solution behavior under varying heating rate and heating load requirements.  First, in Section~\ref{sect:performance-HR} variations due to decreasing heating rate requirements were studied while the heating load remained unconstrained.  It was observed that the peak altitude decreased, the pierce point flight path angle became more shallow, and the total entry time increased as $\dot{Q}_{\max}$ was lowered.  In addition, it was observed that $\dot{Q}_{\max}$ could be lowered to about $6~\textrm{MW}/\textrm{m}^2$ with only slight changes in the performance index value and the angle of attack profile at the start of entry.  Interestingly, only the solutions attained at $\dot{Q}_{\max} = 6~\textrm{MW}/\textrm{m}^2$ and lower saw the peak altitude constrained at its lower limit of $100~\textrm{km}$.  With further reductions in the peak altitude no longer possible, the performance index and the angle of attack near the start of entry began increasing more rapidly as $\dot{Q}_{\max}$ was lowered from $6~\textrm{MW}/\textrm{m}^2$ until convergence failure at $3~\textrm{MW}/\textrm{m}^2$.  The aforementioned trends make sense, because $\dot{Q}$ is a function of both altitude and speed, and the speed is initially high (about $7.35~\textrm{km}/\textrm{s}$) at the start of entry.  Thus, a lower peak altitude, more shallow pierce point flight path angle, and a high angle of attack all contribute to the entry vehicle being able to maintain a sufficiently high altitude during the initial descent into the Earth's atmosphere at high speeds, thereby maintaining $\dot{Q} \leq \dot{Q}_{\max}$.

Somewhat opposite trends were observed in Section~\ref{sect:performance-Q} as the heating load limit was reduced while the heating rate was unconstrained.  As the $Q_{\max}$ limit was lowered, the peak altitude increased, the pierce point flight path angle became more steep, and the total entry time decreased.  However, similar to Section~\ref{sect:performance-HR}, large increases in cost and changes to the angle of attack profile were not observed as $Q_{\max}$ was lowered until about $Q_{\max} = 2600~\textrm{MJ}/\textrm{m}^2$.  Interestingly, $Q_{\max} = 2600~\textrm{MJ}/\textrm{m}^2$ demarcates the transition to solutions where the optimal peak altitude reaches the prescribed upper limit of $200~\textrm{km}$.  Thus, even though opposite bounds were reached, both Sections~\ref{sect:performance-HR} and \ref{sect:performance-Q} suggest that less aggressive entry maneuvers may be obtained by relaxing the bounds placed on the peak altitude.

Despite the somewhat opposite trends observed in Sections~\ref{sect:performance-HR} and \ref{sect:performance-Q}, Section~\ref{sect:performance-combined} demonstrates that significant room exists for reducing $Q_{\max}$ (relative to the natural value of $Q$ obtained when $Q_{\max} = \infty$) even at low values of $\dot{Q}_{\max}$.  Section~\ref{sect:performance-combined} also showed that the trends in peak altitude, pierce point flight path angle, and total entry time observed in Sections~\ref{sect:performance-HR} and \ref{sect:performance-Q} reappeared as $\dot{Q}_{\max}$ was lowered with $Q_{\max}$ held constant, and vice-versa.  For instance, the peak altitude increased and the pierce point flight path angle and total entry time decreased from Mission 1 to Mission 2 as $Q_{\max}$ was lowered from $3500~\textrm{MJ}/\textrm{m}^2$ to $2500~\textrm{MJ}/\textrm{m}^2$ and $\dot{Q}_{\max}$ was held constant at $8~\textrm{MW}/\textrm{m}^2$.  Thus, lessons learned in Sections~\ref{sect:performance-HR} and \ref{sect:performance-Q} may prove useful when designing missions with both heating rate and heating load constraints active.

\section{Conclusions\label{sect:Conclusion}}
The problem of trajectory design for a combined ascent-entry mission with vertical takeoff and impact has been considered.  A trajectory has been designed that includes both boost and entry phases, and Euler parameters were employed in the 3DOF model to parameterize translational motion in the vertical flight phases.  The trajectory design was optimized by numerically solving an optimal control problem using an adaptive Gaussian quadrature collocation method.  It was found that an appropriately chosen performance index could produce trajectories with wide control margins while simultaneously limiting phugoid oscillations during entry.  In addition, it was found that the peak altitude tended to decrease, the pierce point flight path angle became more shallow, and the total entry time increased as the maximum stagnation point heating rate limit was reduced.  In contrast, trajectories obtained at lower values of the maximum stagnation point heating load limit saw the opposite trend.  Overall, as either the heating rate limit or heating load limit was varied, the changes in the trajectory for the boost phases were relatively minor compared to the changes observed in the entry phases.  However, despite changes in the boost phases being relatively minor in comparison, it was found that these changes were critical in enabling a feasible entry trajectory that met both heating rate and heating load constraints.  In addition, adjustments made to the ascent profile enabled, in several instances, trajectories of nearly identical performance even as the heating rate and/or heating load constraints were tightened.

\section*{Acknowledgments}
The authors gratefully acknowledge support for this research from the from the U.S.~National Science Foundation under grants CMMI-1563225, DMS-1522629, and DMS-1819002, from the U.S.~Office of Naval Research under grant N00014-19-1-2543, and from the U.S.~Department of Defense under the National Defense Science \& Engineering Graduate Fellowship (NDSEG) Program.

\renewcommand{\baselinestretch}{1}
\normalsize\normalfont

\bibliographystyle{aiaa}     

\renewcommand{\baselinestretch}{1.5}
\normalsize\normalfont

\end{document}